\magnification=1000
\hsize=11,7cm
\vsize=18.9cm
\lineskip2pt \lineskiplimit2pt
\nopagenumbers

\hoffset=-1truein
\voffset=-1truein

\advance\voffset by 4truecm
\advance\hoffset by 4.5truecm

\newif\ifentete

\headline{\ifentete\ifodd	\count0 
      \rlap{\head}\hfill\tenrm\llap{\the\count0}\relax
    \else
        \tenrm\rlap{\the\count0}\hfill\llap{\head} \relax
    \fi\else
\global\entetetrue\fi}

\def\entete#1{\entetefalse\gdef\head{#1}}
\entete{}

\input amssym.def
\input amssym.tex

\def\-{\hbox{-}}
\def\.{{\cdot}}
\def\O{{\cal O}}
\def\K{{\cal K}}
\def\F{{\cal F}}

\def\Q{{\cal Q}}
\def\G{{\cal G}}

\def\R{{\cal R}}
\def\I{{\cal I}}

\def\ch{\frak c\frak h}

\def\Gr{\frak G\frak r}

\def\int{\frak i\frak n\frak t}

\def\qq{\quad{\rm and}\quad}

\def\mod{\frak m\frak o\frak d}

\def\too{\longrightarrow}
\def\aut{\frak a\frak u\frak t}

 3
 2
\font\large=cmr10  scaled \magstep 2
 2
\font\larti=cmti10  scaled \magstep 2
 2
\font\cds=cmr7
\font\cdt=cmti7

\centerline{\large Equivariant Alperin-Robinson's Conjecture}
\smallskip
\centerline{\large reduces to almost-simple {\larti k}*-groups{\footnote{\dag}{\cds 
We thank Britta Sp\"ath for pointing us a mistake in an earlier version of this paper.}} }
\medskip
\centerline{\bf Lluis Puig }
\medskip
\noindent 
\centerline{\cds CNRS, Institut de Math\'ematiques de Jussieu}
\par
\noindent
\centerline{\cds 6 Av Bizet, 94340 Joinville-le-Pont, France}
\par
\noindent
\centerline{\cds puig@math.jussieu.fr}

\bigskip
\noindent
{\bf Abstract:} {\cds In a recent paper, Gabriel Navarro and Pham Huu Tiep show that the so-called Alperin Weight Conjecture can be verified via the Classification of the Finite Simple Groups, provided any simple group fulfills a very precise list of conditions. Our purpose here is to show that the {\cdt equivariant\/} refinement 
of the Alperin's Conjecture for blocks formulated by Geoffrey Robinson in the eighties can be reduced to checking the {\cdt same\/} statement on any central {\cdt k}*-extension of any finite almost-simple group, or of any finite simple group up to verifying an ``almost necessary'' condition. In an Appendix we develop some old arguments that we need in the proof.}

\bigskip
\noindent
{\bf £1. Introduction }
\medskip
£1.1. In a recent paper [3], Gabriel Navarro and Pham Huu Tiep show that the so-called Alperin Weight Conjecture
can be verified {\it via\/}  the {\it Classification of the Finite Simple Groups\/}, provided any simple group fulfills
a very precise list of conditions that they consider easier to check than ours, firstly stated in [6,~Theorem~16.45] and
significantly weakened in [8,~Theorem~1.6]{\footnote{\dag\dag}{\cds  Gabriel Navarro and Pham Huu Tiep pointed out to us that, when submitting [3], they were not aware of our paper [8], 
only available in arXiv since April 2010.}}. As a matter of fact, our reduction result concerns {\it Alperin's Conjecture for blocks\/} in an {\it equivariant\/} formulation which goes back to Geoffrey Robinson in the eighties (it appears in his joint work [11] with Reiner Staszewski).

\medskip
£1.2. Actually, in the introduction of [6] --- from I29 to I37 --- we consider a refinement of Alperin-Robinson's Conjecture for blocks; but, only in~[8] we really show that its verification can be reduced to check that the {\it same\/} refinement holds on the so-called {\it almost-simple $k^*\-$groups\/} --- namely, central $k^*\-$extensions  of finite  groups $H$
containing a normal simple subgroup $S$ such that $H/S$ is a cyclic $p'\-$group and we have $C_H (S) = \{1\}\,.$   To carry out this checking obviously depends on admitting the  {\it Classification of the Finite Simple Groups\/}, and our proof of the reduction itself uses the {\it solvability\/} of the {\it outer automorphism group\/} of a finite simple group (SOFSG), a known fact whose actual proof depends on this classification.

\medskip
£1.3. Our purpose here is,  from our results in [6] and [8] that we will recall as far as possible, to show that the Alperin-Robinson's Conjecture  for blocks  
can be reduced to checking the {\it same\/} statement on any  almost-simple $k^*\-$group $\hat H$ and moreover, that it may be still reduced to any $k^*\-$central extension of any finite simple group
provided we check an ``almost necessary''   condition (see Proposition~£2.14 below) in such an~$\hat H\,.$ We add an Appendix which actually deals with a more general
situation, but provides tools for the proof of our reduction.

\medskip
£1.4. Explicitly, let $p$ be a prime number, $k$ an algebraically closed  field of characteristic $p\,,$
 $\O$ a complete discrete valuation ring of characteristic zero admitting $k$ as the {\it residue\/} field, and 
$\K$ the field of fractions of~$\O\,.$ Moreover, let $\hat G$ be a $k^*\-$central extension of a finite group~$G$
--- simply called {\it finite $k^*\-$group\/} of $k^*\-$quotient $G$ [6,~1.23] --- $b$ a block of $\hat G$ 
[6,~1.25] and $\G_k (\hat G,b)$  the {\it scalar extension\/} from $\Bbb Z$ to $\O$ of the {\it Grothendieck group\/} of the category of finitely generated $k_*\hat Gb\-$modules [6,~14.3].

\medskip
£1.5. Choose  a maximal Brauer $(b,\hat G)\-$pair $(P,e)\,;$ denote by $\F_{\!(b,\hat G)}$ the category
--- called the {\it Frobenius $P\-$category\/} of $(b,\hat G)$ [6,~3.2]  --- formed by
all the subgroups of $P$ and, if $Q$ and $R$ are subgroups of~$P\,,$ by the group homomorphisms 
$\F_{\!(b,\hat G)}(Q,R)$ from $R$ to $Q$ determined by all the elements $x\in G$ fulfilling $(R,g)\i (Q,f)^x$ where 
$(Q,f)$ and $(R,g)$ are the corresponding  Brauer $(b,\hat G)\-$pairs contained in $(P,e)\,;$ in particular, we set
$$\F_{\!(b,\hat G)}(Q) =\F_{\!(b,\hat G)}(Q,Q)\cong N_G (Q,f)/C_G (Q)
\eqno £1.5.1.$$
 Recall that the Brauer
$(b,\hat G)\-$pair  $(Q,f)$ is called {\it selfcentralizing\/} if the image $\bar f$ of $f$ in $\bar C_{\hat G}(Q) =C_{\hat G}(Q)/Z(Q)$ is a block of  {\it defect zero\/} and then we denote by $(\F_{\!(b,\hat G)})^{^{\rm sc}}$  the full subcategory  of~$\F_{\!(b,\hat G)}$ over the {\it selfcentralizing\/} Brauer $(b,\hat G)\-$pairs contained in $(P,e)\,.$

\medskip
£1.6. Recall that an {\it $(\F_{\!(b,\hat G)})^{^{\rm sc}}\-$chain\/} is just a functor 
$\frak q\,\colon \Delta_n\to (\F_{\!(b,\hat G)})^{^{\rm sc}}$ from the $n\-$simplex $\Delta$ considered as a
category with the morphisms given by the order relation; then, the  {\it proper category of 
$(\F_{\!(b,\hat G)})^{^{\rm sc}}\-$chains\/} --- denoted by $\ch^*\big((\F_{\!(b,\hat G)})^{^{\rm sc}}\big)$
---  is formed by the {\it $(\F_{\!(b,\hat G)})^{^{\rm sc}}\-$chains\/} as objects and by the pairs of order-preserving
maps and natural isomorphisms of functors as morphisms  [6,~A2.8]. Denoting by $\Gr$ the category of finite groups, 
we actually have a functor [6,~Proposition~A2.10]
$$\aut_{(\F_{\!(b,\hat G)})^{^{\rm sc}}} : \ch^*\big((\F_{\!(b,\hat G)})^{^{\rm sc}}\big)\too \Gr
\eqno £1.6.1\phantom{.}$$
mapping the  $(\F_{\!(b,\hat G)})^{^{\rm sc}}\-$chain $\frak q$ on the stabilizer $\F_{\!(b,\hat G)}(\frak q)$
of $\frak q$ in $\F_{\!(b,\hat G)} \big(\frak q (n)\big)\,.$
Moreover, setting $Q = \frak q (n)$ and denoting by $f$ the block of $C_{\hat G}(Q)$ such that
$(P,e)$ contains $(Q,f)\,,$ it is clear that $N_{\hat G}(Q,f)$ acts on the simple $k\-$algebra 
$k_*\bar C_{\hat G}(Q)\bar f$ and it is well-known that this action determines a central $k^*\-$extension
$\hat \F_{\!(b,\hat G)} (Q)$ of $\F_{\!(b,\hat G)}(Q)$ [6,~7.4]; in particular, by restriction we get a central
$k^*\-$extension $\hat\F_{\!(b,\hat G)}(\frak q)$ of $\F_{\!(b,\hat G)}(\frak q)\,.$

\medskip
£1.7.  Denoting by $k^*\-\Gr$
the category of $k^*\-$groups with finite $k^*\-$quotient, in [6,~Theorem~11.32] we prove the existence
of a suitable  functor
$$\widehat\aut_{(\F_{\!(b,\hat G)})^{^{\rm sc}}} : \ch^*\big((\F_{\!(b,\hat G)})^{^{\rm sc}}\big)\too k^*\-\Gr
\eqno £1.7.1\phantom{.}$$
lifting $\aut_{(\F_{\!(b,\hat G)})^{^{\rm sc}}}$ and mapping $\frak q$ on $\hat\F_{\!(b,\hat G)}(\frak q)\,;$
then, still denoting by $\G_k$ the functor mapping any $k^*\-$group with finite $k^*\-$quotient $\hat G$
on   the {\it scalar extension\/} from $\Bbb Z$ to $\O$ of the {\it Grothendieck group\/} of the category of finitely generated $k_*\hat G\-$modules, and any $k^*\-$group homomorphism on the corresponding restriction,
we consider the inverse limit
$$\G_k (\F_{\!(b,\hat G)},\widehat\aut_{(\F_{\!(b,\hat G)})^{^{\rm sc}}}) = \lim_{\longleftarrow}\,(\G_k\circ 
\widehat\aut_{(\F_{\!(b,\hat G)})^{^{\rm sc}}})
\eqno £1.7.2,$$
 called the {\it Grothendieck group of\/} $\F_{\!(b,\hat G)}$ [6,~14.3.3 and~Corollary~14.7]; it follows from [6,~I32 and Corollary~14.32] that Alperin's Conjecture for blocks is actually equivalent to the existence of an $\O\-$module isomorphism  $$\G_k (\hat G,b)\cong \G_k (\F_{\!(b,\hat G)},\widehat\aut_{(\F_{\!(b,\hat G)})^{^{\rm sc}}})
 \eqno £1.7.3\phantom{.}$$
 which actually amounts to saying that both members have the same $\O\-$rank.

 \medskip
 £1.8. Denote by ${\rm Out}_{k^*}(\hat G)$ the group of {\it outer\/} $k^*\-$automorphisms of $\hat G$ and by 
${\rm Out}_{k^*}(\hat G)_b$ the stabilizer of $b$ in ${\rm Out}_{k^*}(\hat G)\,;$ on the one hand, it is clear that 
${\rm Out}_{k^*}(\hat G)_b$ acts on $\G_k (\hat G,b)\,;$ on the other hand, an easy {\it Frattini argument\/} 
shows that the stabilizer ${\rm Aut}_{k^*}(\hat G)_{(P,e)}$ of $(P,e)$ in ${\rm Aut}_{k^*}(\hat G)_b$
{\it covers\/} ${\rm Out}_{k^*}(\hat G)_b$ and it is clear that it acts on $(\F_{\!(b,\hat G)})^{^{\rm sc}}\,,$
so that finally ${\rm Out}_{k^*}(\hat G)_b$  still acts on the inverse limit $\G_k (\F_{\!(b,\hat G)},\widehat\aut_{(\F_{\!(b,\hat G)})^{^{\rm nc}}})$ [6,~16.3 and~16.4]. A stronger question is whether or not in~£1.7.3 there exists a 
{\it ${\rm Out}_{k^*}(\hat G)_b\-$stable\/} isomorphism and in [8,~Theorem~1.6] we prove that it suffices to check 
this statement in the {\it almost-simple $k^*\-$groups\/} considered above.

\medskip
£1.9. Here, we are interested in a weaker form of this question, namely in whether or not
there exists a $\K {\rm Out}_{k^*}(\hat G)_b\-$module isomorphism 
$$\K\otimes_\O\G_k (\hat G,b)\cong \K\otimes_\O\G_k (\F_{\!(b,\hat G)},
\widehat\aut_{(\F_{\!(b,\hat G)})^{^{\rm sc}}})
 \eqno £1.9.1;$$
as a matter of fact, it is a {\it numerical\/} question since it amounts to saying that the 
${\rm Out}_{k^*}(\hat G)_b\-$characters of both members coincide and note that they are actually {\it rational\/} characters. Thus, it makes sense to relate this statement with the old Robinson's {\it equivariant condition\/} 
recalled below. We still need some notation; for any Brauer $(b,\hat G)\-$pair $(Q,f)$ contained in $(P,e)\,,$ 
the group $\F_Q (Q)$ of inner automorphisms of $Q$ is a normal subgroup of $\F_{\!(b,\hat G)}(Q)$
and we set (cf.~£1.5.1)
$$\tilde\F_{\!(b,\hat G)}(Q) = \F_{\!(b,\hat G)}(Q)/\F_Q (Q)\cong N_G (Q,f)/Q\.C_G (Q)
\eqno £1.9.2;$$
\eject
\noindent
moreover, if $(Q,f)$ is selfcentralizing then $\F_Q (Q)$ can be identified to a normal $p\-$subgroup of 
$\hat\F_{\!(b,\hat G)}(Q)\,;$ then, we also set
$$\skew4\hat{\tilde\F}_{\!(b,\hat G)} (Q) = \hat\F_{\!(b,\hat G)}(Q)/\F_Q (Q)
\eqno £1.9.3\phantom{.}$$
and  denote by $o_{(Q,f)}$  the sum of blocks of {\it defect zero\/} of $\skew4\hat{\tilde\F}_{\!(b,\hat G)} (Q)\,;$
note that, since  the stabilizer ${\rm Aut}_{k^*}(\hat G)_{(P,e)}$ of $(P,e)$ in ${\rm Aut}_{k^*}(\hat G)_b$
{\it covers\/} ${\rm Out}_{k^*}(\hat G)_b$ and  acts on $(\F_{\!(b,\hat G)})^{^{\rm sc}}\,,$
 the stabilizer  $C_{(Q,f)}$ in a (cyclic) subgroup $C$ of ${\rm Out}_{k^*}(\hat G)_b$ of the $G\-$conjugacy 
 class of~$(Q,f)$ acts naturally on $\G_k \big(\skew4\hat{\tilde\F}_{\!(b,\hat G)}  (Q),o_{(Q,f)}\big)\,.$

\medskip
£1.10. Following Robinson, let us consider the following {\it equivariant condition\/}:
\smallskip
\noindent
(E)\quad {\it For any cyclic subgroup $C$ of ${\rm Out}_{k^*}(\hat G)_b$ we have
$${\rm rank}_\O\big(\G_k (\hat G,b)^C\big) = \sum_{(Q,f)} {\rm rank}_{\O} 
\Big(\G_k \big(\skew4\hat{\tilde\F}_{\!(b,\hat G)}  (Q),o_{(Q,f)}\big)^{C_{(Q,f)}}\Big)
\eqno £1.10.1\phantom{.}$$
where $(Q,f)$ runs over a set of representatives contained in $(P,e)$ for the set of $C\-$orbits of 
$G\-$conjugacy classes of selfcentralizing  Brauer $(b,\hat G)\-$pairs  and, for such a $(Q,f)\,,$ 
we denote by $C_{(Q,f)}$ the stabilizer in~$C$ of the $G\-$conjugacy class of~$(Q,f)\,.$\/}
\smallskip
\noindent
We are ready to state our first main result.

\bigskip 
\noindent
{\bf Theorem~£1.11.} {\it Assume {\rm (SOSFG)} and that any block  $c$ of any $k^*\-$extension $\hat H$ 
 of any finite  group $H\,,$ containing a finite simple group $S$ such that $H/S$ is a cyclic 
 $p'\-$group and that we have  $C_H (S) = \{1\}\,,$  fulfills the equiva-riant condition {\rm (E)}. Then,  any block  $b$ of any $k^*\-$extension $\hat G$ 
of any finite group $G$  fulfills the equivariant condition {\rm (E)} and, in particular, we have a
 $\K {\rm Out}_{k^*}(\hat G)_b\-$module isomorphism 
$$\K\otimes_\O \G_k (\hat G,b)\cong \K\otimes_\O 
\G_k (\F_{\!(b,\hat G)},\widehat\aut_{(\F_{\!(b,\hat G)})^{^{\rm sc}}})
\eqno £1.11.1.$$\/}

\bigskip
\noindent
{\bf £2. The obstruction}

\bigskip
£2.1. In order to explain the obstruction to get a better reduction, let $\hat S$ be a  $k^*\-$group
of non-abelian simple $k^*\-$quotient $S$ and $d$ a block of $\hat S$ which fulfill condition (E);
choose a maximal Brauer $(d,\hat S)\-$pair $(P,e)$ and denote by~$\Q$ a set of representatives contained in 
$(P,e)$ for the set of $S\-$conjugacy classes of selfcentralizing Brauer $(d,\hat S)\-$pairs;
   once again  since  the stabilizer 
 ${\rm Aut}_{k^*}(\hat S)_{(P,e)}$ of $(P,e)$ in ${\rm Aut}_{k^*}(\hat S)_d$
{\it covers\/} ${\rm Out}_{k^*}(\hat S)_d$ and  acts on $(\F_{\!(d,\hat S)})^{^{\rm sc}}\,,$ the group
${\rm Out}_{k^*}(\hat S)_d$ acts on the family $\{\G_k \big(\skew4\hat{\tilde\F}_{\!(d,\hat S)} 
 (Q),o_{(Q,f)}\big)\}_{(Q,f)\in \Q}$ and\break
 \eject
 \noindent
  the direct sum of this family becomes an 
 $\O {\rm Out}_{k^*}(\hat S)_d\-$module. Then, since both $\K{\rm Out}_{k^*}(\hat S)_d\-$modules
$$\K\otimes_\O \G_k (\hat S,d)\qq \bigoplus_{(Q,f)\in \Q} \K\otimes_\O
\G_k \big(\skew4\hat{\tilde\F}_{\!(d,\hat S)}  (Q),o_{(Q,f)}\big)
\eqno £2.1.1\phantom{.}$$actually come from $\Bbb Q{\rm Out}_{k^*}(\hat S)_d\-$modules,  equalities~£1.10.1 amount to saying that 
these ${\rm Out}_{k^*}(\hat S)_d\-$representations have the same character and therefore that we have
a $\K{\rm Out}_{k^*}(\hat S)_d\-$module isomorphism
$$\K\otimes_\O \G_k (\hat S,d)\cong \bigoplus_{(Q,f)\in \Q} \K\otimes_\O\G_k 
\big(\skew4\hat{\tilde\F}_{\!(d,\hat S)}  (Q),o_{(Q,f))}\big)
\eqno £2.1.2.$$

\medskip
£2.2. Let $\hat H$ be a $k^*\-$group of finite $k^*\-$quotient $H$ in such a way that $\hat H$ contains
and normalizes $\hat S\,,$ and $c$ a block of $\hat H$ such that $cd \not= 0\,;$ denoting by $\hat H_d$
the stabilizer of $d$ in $\hat H\,,$ {\it Fong's reduction\/} can be written as follows [7,~Propositions~3.2 and~3.5]
$$k_*\hat Hc = {\rm Ind}_{\hat H_d}^{\hat H} \big(k_*\hat H_d(cd)\big)
\eqno £2.2.1\phantom{.}$$
and we know that $cd$ is a block of $\hat H_d\,;$ hence, for our purposes, we may assume that $\hat H$ fixes~$d$ and thus that we actually have $cd = c\,.$  As above, we assume that
$A = H/S$ is a cyclic $p'\-$group and that we have $C_H (S) = \{1\}\,.$

\medskip
£2.3.  On the other hand,  it is clear that $C_{\hat H}(P,e)$ acts on the $k^*\-$group $\skew4\hat{\tilde\F}_{\!(d,\hat S)}  (P)$ acting trivially on ${\tilde\F}_{\!(d,\hat S)}  (P)\,;$ let us denote by 
$K_{\hat H}(P,e)$ the kernel of this action and set
$$\hat L = \hat S\.K_{\hat H}(P,e)\qq D =\hat L/\hat S
\eqno £2.3.1;$$
it follows from [8,~Proposition~3.8] that $c$ is still a block of $\hat L$ and from [8,~Theorem~3.10] that
the {\it source $P\-$interior algebras\/} of $(d,\hat S)$ and $(c,\hat L)$ are isomorphic; in particular, we have 
[4,~Propositions~6.12 and~14.6]
$$\skew4\hat{\tilde\F}_{\!(d,\hat S)}  (P) = \skew4\hat{\tilde\F}_{\!(c,\hat L)}  (P)
\eqno £2.3.2.$$

\medskip
£2.4. More precisely, it follows from [6,~Lemma~15.16] that the block $e$ of $C_{\hat S}(P)$ splits into a family
$\{e_\varphi\}_{\varphi\in {\rm Hom}(D,k^*)}$ of blocks of 
$$C_{\hat L}(P,e) = K_{\hat H}(P,e)
\eqno £2.4.1\phantom{.}$$ 
and then any $(P,e_\varphi)$ 
clearly becomes a maximal Brauer $(d_\varphi,\hat L)\-$pair for a suitable block $d_\varphi$ of $\hat L\,;$
since, by the very definition of $K_{\hat H}(P,e)\,,$ the group ${\tilde\F}_{\!(d,\hat S)}  (P)$ acts trivially
on this $k^*\-$group, the number of blocks $d_\varphi$ coincides with $\vert {\rm Hom}(D,k^*)\vert$
and then, a simple argument on the dimensions proves that
$$k_*\hat S d\cong k_*\hat L d_\varphi
\eqno £2.4.2\phantom{.}$$
\eject
\noindent
for any $\varphi\in {\rm Hom}(D,k^*)\,.$ Moreover, from £2.3 above, we have $d_\varphi = c$
for some choice of $\varphi\,;$ set $e_c = e_\varphi$ and note that it follows from equality~£2.3.2 above and from
[8,~Corollary~3.12] that we have
$$ \F_{\!(d,\hat S)}   = \F_{\!(c,\hat L)} \qq \widehat\aut_{(\F_{\!(d,\hat S)})^{^{\rm sc}}} 
= \widehat\aut_{(\F_{\!(c,\hat L)})^{^{\rm sc}}}
\eqno £2.4.3.$$
Hence, for our purposes,  we may replace the pair $(d,\hat S)$  by the pair $(c,\hat L)\,;$ in particular, 
denoting by $\R$ a set of representatives contained in $(P,e_c)$ for the set 
of $L\-$conjugacy classes of selfcentralizing Brauer $(c,\hat L)\-$pairs,
from isomorphism~£2.1.2 we get a $\K{\rm Out}_{k^*}(\hat L)_c\-$module isomorphism
$$\K\otimes_\O \G_k (\hat L,c)\cong \bigoplus_{(R,g)\in \R} \K\otimes_\O\G_k 
\big(\skew4\hat{\tilde\F}_{\!(c,\hat L)}  (R),o_{(R,g)}\big)
\eqno £2.4.4\phantom{.}$$
since ${\rm Aut}_{k^*}(\hat L)$ clearly stabilizes $\hat S\,.$

\medskip
£2.5.  Set $B = N_{\hat H}(P,e_c)/N_{\hat L}(P,e_c)\cong \hat H/\hat L\,;$ isomorphism~£2.4.4 is obviously a $\K B\-$isomorphism and therefore, since the group 
$B$ is cyclic, the respective $B\-$stable $\K\-$bases
$${\rm Irr}_k (\hat L,c)\qq \bigsqcup_{(R,g)\in \R} {\rm Irr}_k 
\big(\skew4\hat{\tilde\F}_{\!(c,\hat L)}  (R),o_{(R,g)}\big)
\eqno £2.5.1\phantom{.}$$
become {\it isomorphic $B\-$sets\/}. That is to say, an irreducible Brauer character $\theta$ of $\hat L$
 in the block $c$ determines a  selfcentralizing  Brauer $(c,\hat L)\-$pair $(R,g)$ in~$\R$ and
a {\it projective\/} irreducible Brauer character $\theta^*$ of $\skew4\hat{\tilde\F}_{\!(c,\hat L)}  (R)\,,$
in such a way that the stabilizer $B_\theta$ of $\theta$ in $B$ coincides with the stabilizer of the pair formed 
by~$(R,g)\in \R$ and $\theta^*\,.$

\medskip
£2.6. On the one hand, note that $c$ is also a block of the stabilizer $\hat H_\theta$ of $\theta\in {\rm Irr}_k (c,\hat L)$  in $\hat H\,;$ denote by  $\G_k (\hat H_\theta\!\mid\! \theta)$ the direct summand of 
$\G_k (\hat H_\theta,c)$ generated by the classes of the simple $k_*\hat H_\theta\-$modules whose restriction to
$\hat L$ involves $\theta\,;$ then, it follows from the so-called {\it Clifford theory\/} that we have a canonical isomorphism
$$\G_k (\hat H_\theta\!\mid\! \theta)\cong \G_k (\widehat B^{^\theta}_\theta)
\eqno £2.6.1\phantom{.}$$
for the central $k^*\-$extension $\widehat B^{^\theta}_\theta$ of $B_\theta$ defined in~£2.7.3 below. Moreover, recall that we have
$$\G_k (\hat H,c) = \bigoplus_{\theta\in \Theta} {\rm Ind}_{\hat H_\theta}^{\hat H} 
\big(\G_k (\hat H_\theta\!\mid\! \theta)\big)
\eqno £2.6.2\phantom{.}$$
where $\Theta$ is a set of representatives for the set of $H\-$orbits of ${\rm Irr}_k (\hat L,c)\,.$
 Consequently, since any cyclic subgroup $C$ of ${\rm Out}_{k^*}(\hat H)_c$  acts on ${\rm Irr}_k(\hat L)\,,$  we have
 $${\rm rank}_\O\big(\G_k (\hat H,c)^C\big) = \sum_{\theta\in \Theta} 
 {\rm rank}_\O \big(\G_k (\widehat B^{^\theta}_\theta)^{C_\theta}\big)
\eqno £2.6.3\phantom{.}$$
where, for any $\theta\in \Theta\,,$ $C_\theta$ denotes the stabilizer in~$C$ of~ the $H\-$orbit of $\theta\,.$
\eject

\medskip
£2.7.  More explicitly, denoting by 
$V_\theta$ a $k_*\hat L\-$module affording $\theta$ and by $\rho_\theta\,\colon \hat L\to GL_k (V_\theta)$
the corresponding $k^*\-$group homomorphism, the action of $H_\theta$ on $\hat L$ determines
a group homomorphism $\bar\rho_\theta\,\colon H_\theta\to PGL_k (V_\theta)$ and we can consider the {\it pull-back\/}
$$\matrix{H_\theta&\buildrel \bar\rho_\theta\over\too &PGL_k (V_\theta)\cr
\uparrow&&\uparrow{\scriptstyle \pi_\theta}\cr
\widehat H_\theta^{^\theta}&\too &GL_k (V_\theta)}
\eqno £2.7.1\phantom{.}$$
where $\widehat H_\theta^{^\theta}$ is the $k^*\-$group formed by the pairs $(x,f)\in H_\theta\times GL_k (V_\theta)$
such that $\bar\rho_\theta (x) = \pi_\theta (f)\,.$
 Moreover, since the composition of $\bar\rho_\theta$ with the map $\hat L\to H_\theta$
 extends $\pi_\theta\circ\rho_\theta\,,$ we have an injective canonical $k^*\-$group homomorphism $\hat L\to \widehat H_\theta^{^\theta}\,,$
so that we get an injective canonical $k^*\-$group homomorphism
$$L\times k^*\cong \hat L * (\hat L)^\circ \too \hat H_\theta * (\widehat H_\theta^{^\theta})^\circ
\eqno £2.7.2\phantom{.}$$
 and, identifying $L$ with its image which is a normal subgroup of $\hat H_\theta * (\widehat H_\theta^{^\theta})^\circ\,,$ we set
$$\widehat B^{^\theta}_\theta = \big(\hat H_\theta * (\widehat H_\theta^{^\theta})^\circ\big)\big/ L
\eqno £2.7.3.$$

\medskip
£2.8. On the other hand, for any selfcentralizing  Brauer  $(c,\hat L)\-$pair $(R,g)$ contained in~$(P,e_c)\,,$ 
set $\bar N_{\hat L} (R,g) = N_{\hat L} (R,g)/R$ 
and denote by $\bar g$ the image of $g$ 
in $k_*\bar N_{\hat L} (R,g)\,;$ since $\bar g$ is a block of {\it defect zero\/} of $\bar C_{\hat L}(R)\,,$
applying again {\it Fong's reduction\/} we get [7,~Proposition~3.2 and~Theorem~3.7]
$$k_*\bar N_{\hat L} (R,g)\bar g \cong k_*\bar C_{\hat L} (R)\bar g\otimes_k 
k_* \skew4\hat{\tilde\F}_{\!(c,\hat L)}  (R)
\eqno £2.8.1\phantom{.}$$
and clearly $\bar g\otimes o_{(R,g)}$ corresponds to the sum of all the blocks  of {\it defect zero\/} 
in~$k_*\bar N_{\hat L} (R,g)\bar g\,.$ Note that, since the quotient $\bar N_{\hat H}(R,g)/\bar N_{\hat L}(R,g)$
is a $p'\-$group, the blocks of {\it defect zero\/} of $\bar N_{\hat H}(R,g)$ and $\bar N_{\hat L}(R,g)$ mutually
correspond [6,~Proposition~15.9].

\medskip
£2.9. Moreover, since the quotient $E = C_{\hat H}(R,g)/C_{\hat L}(R)$ is cyclic,
it follows again from [6,~Lemma~15.16] that the block $g$ of $C_{\hat L}(R)$ splits into a family
$\{g_\psi\}_{\psi\in {\rm Hom}(E,k^*)}$ of blocks of $C_{\hat H}(R,g)$ and, according to
[8,~3.7], the group ${\tilde\F}_{\!(c,\hat L)}  (R)$ acts transitively on this family; as in~£2.4 above, a simple 
argument on the dimensions proves that
$$k_*C_{\hat L}(R)g \cong k_*C_{\hat H}(R,g)g_\psi
\eqno £2.9.1\phantom{.}$$
for any $\psi\in {\rm Hom}(E,k^*)\,;$ setting $\bar N_{\hat H} (R,g) = N_{\hat H} (R,g)/R\,,$ 
once again  {\it Fong's reduction\/} provides the following 
decomposition [7,~Proposition~3.2 and~Theorem~3.7]
$$k_*\bar N_{\hat H} (R,g)\bar g \cong {\rm Ind}_{\bar N_{\hat H} (R,g_\psi)}^{\bar N_{\hat H} (R,g)}\big(k_*\bar C_{\hat H}(R,g)\bar g_\psi \otimes_k k_* \skew4\hat{\tilde\F}_{\!(c,\hat H)}  (R)\big)
\eqno £2.9.2\phantom{.}$$
where, as in~£1.9.2 above, we have
$${\tilde\F}_{\!(c,\hat H)}  (R)\cong \bar N_{\hat H} (R,g_\psi)/\bar C_{\hat H} (R)
\eqno £2.9.3.$$
\eject

\medskip
£2.10. Furthermore, since ${\rm Hom}(E,k^*)$ is a $p'\-$group, $p$ does not divide 
$\vert \bar N_{\hat H} (R,g)\colon \bar N_{\hat H} (R,g_\psi)\vert$ and therefore isomorphism~£2.9.2 induces a bijection between the sets of isomorphism classes of {\it projective simple\/} 
$k_*\bar N_{\hat H} (R,g)\bar g\-$ and~$k_* \skew4\hat{\tilde\F}_{\!(c,\hat H)}  (R)\-$modules; 
in other words, ${\rm Tr}_{\bar N_{\hat H} (R,g_\psi)}^{\bar N_{\hat H} (R,g)}
(\bar g_\psi\otimes o_{(R,g_\psi)})$  corres-ponds to the sum $n_{(R,g)}$ of all the blocks  of {\it defect zero\/} 
of~$k_*\bar N_{\hat H} (R,g)\bar g\,.$ But, a maximal Brauer  $(c,\hat H)\-$pair $(P,e_\circ)$ 
such that $e_\circ$ appears in the decomposition of $e_c$ (cf.~£2.4) contains $(R,g_\psi)$ for a unique choice of~$\psi\in {\rm Hom}(E_\theta,k^*)\,,$ and we set $g_\psi = g_\circ\,;$ in particular, the family 
$\R_\circ = \{(R,g_\circ)\}_{(R,g)\in \R}$ is  a set of representatives contained in $(P,e_\circ)$ for the set 
of $H\-$conjugacy classes of selfcentralizing Brauer $(c,\hat H)\-$pairs. In conclusion, isomorphism~£2.9.2 induces 
a natural $\O\-$isomorphism
$$\G_k \big(\bar N_{\hat H} (R,g),n_{(R,g)}\big)\cong 
\G_k \big(\skew4\hat{\tilde\F}_{\!(c,\hat H)}  (R), o_{(R,g_\circ)}\big)
\eqno £2.10.1.$$

\medskip
£2.11. Now, denote by $\Theta_{(R,g)}$ the subset of elements of $\Theta$ (cf.~£2.6) determining the Brauer  
$(c,\hat L)\-$pair $(R,g)\,;$ any $\theta\in \Theta_{(R,g)}$ also determines~a {\it projective\/} irreducible
character $\theta^*$ of $\skew4\hat{\tilde\F}_{\!(c,\hat L)}  (R)$ and then, according to isomorphism~£2.8.1,
$\theta^*$ determines a  {\it projective\/} irreducible Brauer character $\zeta_{\theta^*}$ 
of~$k_*\bar N_{\hat L} (R,g)\bar g \,.$ Since $\hat H_\theta$ is also the stabilizer in $\hat H$ of the pair formed by
$(R,g)\in \R$ and $\theta^*$ (cf.~£2.5), we have
$$B_\theta\cong N_{\hat H_\theta}(R,g)/N_{\hat L}(R,g)
\eqno £2.11.1\phantom{.}$$
and, denoting by  $\G_k \big(\bar N_{\hat H_\theta}  (R,g)\!\mid\! \zeta_{\theta^*}\big)$ the direct summand of  
$\G_k \big(\bar N_{\hat H_\theta} (R,g),\bar g\big)$ generated by the classes of simple $k_*\bar N_{\hat H_\theta}  (R,g)\-$modules whose restriction to $\bar N_{\hat L} (R,g)$ involves $\zeta_{\theta^*}\,,$
 it follows 
again from the  {\it Clifford theory\/} that we have a canonical isomorphism
$$\G_k \big(\bar N_{\hat H_\theta}  (R,g)\!\mid\! \zeta_{\theta^*}\big) \cong 
\G_k  (\widehat B^{^{\theta^*}}_\theta)
\eqno £2.11.2\phantom{.}$$
for an analogous central $k^*\-$extension $\widehat B^{^{\theta^*}}_\theta$ of $B_\theta\,.$
 Moreover, since  the blocks of {\it defect zero\/} of $\bar N_{\hat H}(R,g)$ and of $\bar N_{\hat L}(R,g)$ mutually
correspond, as in~£2.6 above we have
$$\G_k \big(\bar N_{\hat H} (R,g),n_{(R,g)}\big) = \bigoplus_{\theta\in \Theta_{(R,g)}} {\rm Ind}_{N_{\hat H_\theta}  (R,g)}^{N_{\hat H}  (R,g)} \Big(\G_k \big(\bar N_{\hat H_\theta}  (R,g)\!\mid\! \zeta_{\theta^*}\big)\Big)
\eqno £2.11.3.$$

\medskip
£2.12.  Finally, it is clear that ${\rm Aut}_{k^*}(\hat H)$ stabilizes $\hat S$ and, consequently, that 
${\rm Aut}_{k^*}(\hat H)_c$ stabilizes $\hat L$  and acts on $\G_k(\hat L,c)\,;$ more precisely,
 the stabilizer ${\rm Aut}_{k^*}(\hat H)_{(P,e_c)}$ of $(P,e_c)$ in  ${\rm Aut}_{k^*}(\hat H)_c$
 acts on $(\F_{\!(c,\hat L)})^{^{\rm sc}}$ and the $\K{\rm Out}_{k^*}(\hat L)_c\-$module isomorphism~£2.4.4 
restricts to a $\K\big({\rm Aut}_{k^*}(\hat H)_{(P,e_c)}\big)\-$iso-morphism; moreover, this isomorphism
clearly preserves the {\it $\K B\-$isotypic components\/} of both members and, since 
${\rm Aut}_{k^*}(\hat H)_{(P,e_c)}$ {\it covers\/} ${\rm Out}_{k^*}(\hat H)_c\,,$
we still have a $\K {\rm Out}_{k^*}(\hat H)_c\-$module isomorphism
$$\big(\K\otimes_\O \G_k (\hat L,c)\big)^B\cong \Big(\bigoplus_{(R,g)\in \R} \K\otimes_\O\G_k 
\big(\skew4\hat{\tilde\F}_{\!(c,\hat L)}  (R),o_{(R,g)}\big)\Big)^B
\eqno £2.12.1;$$
in particular, as in~2.4 above, for any cyclic subgroup $C$ of~${\rm Out}_{k^*} (\hat H)_c\,,$  
the respective $C\-$stable $\K\-$bases indexed by the quotient sets
$${\rm Irr}_k (\hat L,c)/B\qq \Big(\bigsqcup_{(R,g)\in \R} {\rm Irr}_k 
\big(\skew4\hat{\tilde\F}_{\!(c,\hat L)}  (R),o_{(R,g)}\big)\Big)\Big/B
\eqno £2.12.2\phantom{.}$$
become {\it isomorphic $C\-$sets\/}.

\medskip
£2.13. Consequently, for any cyclic subgroup $C$ of~${\rm Out}_{k^*} (\hat H)_c\,,$  we can choose 
a {\it $B\-$set isomorphism\/}
$${\rm Irr}_k (\hat L,c)\cong  \bigsqcup_{(R,g)\in \R} {\rm Irr}_k 
\big(\skew4\hat{\tilde\F}_{\!(c,\hat L)}  (R),o_{(R,g)}\big)
\eqno £2.13.1\phantom{.}$$
inducing a {\it $C\-$set isomorphism\/} between the quotient sets in~£2.12.2; in this situation, 
considering $\theta\in \Theta$ and the corresponding pair formed by $(R,g)\in \R$ and by a {\it projective\/} irreducible
character $\theta^*$ of $\skew4\hat{\tilde\F}_{\!(c,\hat L)}  (R)\,,$ the stabilizer in $C_{(R,g)}$ of $\theta^*$
coincides with $C_\theta$ and we have (cf.~£2.10.1 and~£2.11.2)
$$\eqalign{\G_k \big(\skew4\hat{\tilde\F}_{\!(c,\hat H)}  (R), o_{(R,g_\circ)}\big)^{C_{(R,g_\circ)}}
&\cong \G_k \big(\bar N_{\hat H} (R,g),n_{(R,g)}\big)^{C_{(R,g)}}\cr
&\cong \bigoplus_{\theta\in \Theta_{(R,g)}} \G_k \big(\bar N_{\hat H_\theta}  (R,g)
\!\mid\! \zeta_{\theta^*}\big)^{C_\theta}\cr
&\cong \bigoplus_{\theta\in \Theta_{(R,g)}} \G_k  (\widehat B^{^{\theta^*}}_\theta)^{C_\theta}\cr}
\eqno £2.13.2.$$
 At this point, it follows from equalities~£2.6.3 and~isomorphisms~£2.13.2 that a {\it sufficient\/} statement
 to guaranteeing that the block $c$ of $\hat H$ fulfillis condition (E) is that, for any cyclic subgroup $C$  
 of~${\rm Out}_{k^*} (\hat H)_c$ and any $\theta\in \Theta\,,$ the following equality holds
$$ {\rm rank}_\O \big(\G_k (\widehat B^{^\theta}_\theta)^{C_\theta}\big) = {\rm rank}_\O 
\big(\G_k  (\widehat B^{^{\theta^*}}_\theta)^{C_\theta}\big)
\eqno £2.13.3.$$
Note that this equality forces that the action of $C_\theta$ is trivial on $\G_k (\widehat B^{^\theta}_\theta)$
if and only if it is trivial on $\G_k  (\widehat B^{^{\theta^*}}_\theta)\,.$ We are ready to state our second main result.

\bigskip
\noindent
{\bf Proposition~£2.14.} {\it With the notation and the hypothesis above, assume that we have  a 
$\K{\rm Out}_{k^*}(\hat L)_c\-$module isomorphism
$$\K\otimes_\O \G_k (\hat L,c)\cong \bigoplus_{(R,g)\in \R} \K\otimes_\O\G_k 
\big(\skew4\hat{\tilde\F}_{\!(c,\hat L)}  (R),o_R\big)
\eqno £2.14.1.$$
Let $C$ be a cyclic subgroup of~${\rm Out}_{k^*} (\hat H)_c$ and $\theta$ and element in ${\rm Irr}_k (\hat L,c)\,.$
 If the actions of $C_\theta$ on $\G_k (\hat H_\theta\!\mid\! \theta)$ and on 
$\G_k \big(\bar N_{\hat H_\theta}  (R,g)\!\mid\! \zeta_{\theta^*}\big)$ have the same kernel then we have 
$$ {\rm rank}_\O \big(\G_k (\widehat B^{^\theta}_\theta)^{C_\theta}\big) = 
{\rm rank}_\O \big(\G_k (\widehat B^{^{\theta^*}}_\theta)^{C_\theta}\big)
\eqno £2.14.2.$$\/}

\bigskip
\noindent
{\bf £3. Proof of the second main result}

\bigskip
£3.1. Let $A$ be a cyclic $p'\-$group and set $\hat A = A\times k^*\,;$ we have  an obvious split exact sequence
$$1\too {\rm Hom}(A,k^*)\too {\rm Aut}_{k^*}(\hat A) \too {\rm Aut}(A)\too 1
\eqno £3.1.1\phantom{.}$$
and a canonical $\O\-$module isomorphism
$$\G_k (\hat A)\cong \O  {\rm Hom}_{k^*}(\hat A,k^*)
\eqno £3.1.2;$$
moreover, the action of ${\rm Hom}(A,k^*)\i {\rm Aut}_{k^*}(\hat A)$ on $\O {\rm Hom}_{k^*}(\hat A,k^*)$ 
through  isomorphism~£3.1.2 is just defined by the ``product'' 
$${\rm Hom}(A,k^*)\times {\rm Hom}_{k^*}(\hat A,k^*)\too {\rm Hom}_{k^*}(\hat A,k^*)
\eqno £3.1.3.$$ 
On the other hand, a subgroup of ${\rm Hom}(A,k^*)$ is the image of ${\rm Hom}(A/D,k^*)$ for 
some subgroup $D$ of $A$ and it is easily checked that the restriction induces an  $\O\-$module isomorphism
$$\G_k (\hat A)^{{\rm Hom}(A/D,k^*)} \cong \G_k (\hat D)
\eqno £3.1.4.$$

\bigskip
\noindent
{\bf Lemma~£3.2.} {\it With the notation above, let $C$ and $C'$ be subgroups of ${\rm Aut}_{k^*}(\hat A)$
having the same order and the same image in ${\rm Aut}(A)\,.$ Then, we have 
$$ {\rm rank}_\O \big(\G_k (\hat A)^{C}\big) = {\rm rank}_\O \big(\G_k (\hat A)^{C'}\big)
\eqno £3.2.1.$$\/}

\par\noindent
{\bf Proof:} We argue by induction on $\vert A\vert\,;$ according to our hypothesis, we have
$$\vert C\cap {\rm Hom}(A,k^*)\vert = \vert C'\cap {\rm Hom}(A,k^*)\vert
\eqno £3.2.2;$$
since ${\rm Hom}(A,k^*)$ is cyclic, we actually get the equality
$$ C\cap {\rm Hom}(A,k^*) = C'\cap {\rm Hom}(A,k^*)
\eqno £3.2.3\phantom{.}$$
\eject
\noindent
and this intersection is the image of ${\rm Hom}(A/D,k^*)$ for some subgroup $D$ of~$A\,;$
moreover, the restriction determines a commutative diagram of short exact sequences
$$\matrix{&1&&1&&1\cr
&\downarrow&&\downarrow&&\downarrow\cr
1\too &{\rm Hom}(A/D,k^*)&\too&X&\too &Y&\too 1\cr
&\downarrow&&\downarrow&&\downarrow\cr
1\too& {\rm Hom}(A,k^*)&\too& {\rm Aut}_{k^*}(\hat A) &\too& {\rm Aut}(A)&\too 1\cr
&\downarrow&&\downarrow&&\downarrow\cr
1\too& {\rm Hom}(D,k^*)&\too& {\rm Aut}_{k^*}(\hat D) &\too& {\rm Aut}(D)&\too 1\cr
&\downarrow&&\downarrow&&\downarrow\cr
&1&&1&&1\cr}
\eqno £3.2.4\phantom{.}$$
where all the horizontal sequences are split since we can choose compatible splittings in the middle and the bottom
 horizontal sequences.
 
 \smallskip
 Then, since $C$ and $C'$ have the same image in ${\rm Aut}(A)\,,$ they also have the same image in
 ${\rm Aut}(D)$ and we actually get $C\cap X = C'\cap X\,,$ so that the images $\tilde C$ and $\tilde C'$
 of $C$ and $C'$ in ${\rm Aut}_{k^*}(\hat D)$ still have the same order; moreover, according to £3.1.4 
 and  £3.2.3, we obtain
 $$\G_k (\hat A)^C \cong \G_k (\hat D)^{\tilde C}\qq \G_k (\hat A)^{C'} \cong \G_k (\hat D)^{\tilde C'}
 \eqno £3.2.5;$$
 thus, if $D\not= A$ then it suffices to apply our induction hypothesis. 
 
 \smallskip
 From now on, we assume that
$$ C\cap {\rm Hom}(A,k^*) = \{1\} = C'\cap {\rm Hom}(A,k^*)
\eqno £3.2.6.$$
 Let us consider the {\it residual\/} Grothendieck group of $\hat A$ [6,~15.22]
 $$\R\G_k (\hat A) = \bigcap_E {\rm Ker}\big({\rm Res}_{\hat E}^{\hat A}\big) 
  \eqno £3.2.7\phantom{.}$$
where $E$ runs over the set of proper subgroups of $A$ and, for such an $E\,,$ we denote by
$${\rm Res}_{\hat E}^{\hat A} : \G_k(\hat A)\too \G_k (\hat E)
\eqno £3.2.8\phantom{.}$$
the restriction map; it is easily checked that we have a canonical isomorphism
 [6,~15.23.4]
 $$\G_k (\hat A) \cong \bigoplus_E \R\G_k (\hat E)
 \eqno £3.2.9\phantom{.}$$
where $E$ runs over the set of subgroups of $A\,.$ In particular, we get
$$\G_k (\hat A)^C \cong \bigoplus_E \R\G_k (\hat E)^C\qq \G_k (\hat A)^{C'} \cong \bigoplus_E \R\G_k (\hat E)^{C'}
\eqno £3.2.10\phantom{.}$$
where $E$ runs over the set of subgroups of $A\,.$
\eject

\smallskip
On the other hand, for any proper subgroup $E$ of $A\,,$  $C$ and $C'$ have as above the same image in 
${\rm Aut}(E)\,;$ moreover, the products of $C$ and $C'$ by the image of ${\rm Hom}(A/E,k^*)$ in 
${\rm Aut}_{k^*}(\hat A)$ have the same order, and their images in ${\rm Aut}_{k^*}(\hat E)$ respectively coincide with the images of $C$ and $C'\,;$ now, the argument above proves that these images also have the same order. Consequently, 
it follows again from our induction hypothesis that we already have
$$ {\rm rank}_\O \big(\G_k (\hat E)^{C}\big) =  {\rm rank}_\O \big(\G_k (\hat E)^{C'}\big)
\eqno £3.2.11.$$
In particular, according to the corresponding isomorphisms~£3.2.10, this forces
$$\sum_F {\rm rank}_\O \big(\R\G_k (\hat F)^{C}\big) = \sum_F {\rm rank}_\O \big(\R\G_k (\hat F)^{C'}\big)
\eqno £3.2.12\phantom{.}$$
where $F$ runs over the set of subgroups of $E\,,$ and, since equalities hold for any proper subgroup $E$ of $A\,,$
we finally obtain
$$ {\rm rank}_\O \big(\R\G_k (\hat E)^{C}\big) =  {\rm rank}_\O \big(\R\G_k (\hat E)^{C'}\big)
\eqno £3.2.13.$$

\smallskip
In conclusion, always according to the  isomorphisms £3.2.10, it remains to prove that
$$ {\rm rank}_\O \big(\R\G_k (\hat A)^{C}\big) =  {\rm rank}_\O \big(\R\G_k (\hat A)^{C'}\big)
\eqno £3.2.14.$$
The equality $\hat A = A\times k^*$ defines an splitting of the exact sequence~£3.1.1
$${\rm Aut}_{k^*}(\hat A)\cong {\rm Hom}(A,k^*)\rtimes {\rm Aut}(A)
\eqno £3.2.15\phantom{.}$$
and $\O\-$module isomorphisms
$$\G_k (\hat A)\cong \G_k (A)\qq \R\G_k (\hat A)\cong \R\G_k (A)
\eqno £3.2.16;$$
we make the corresponding identifications; note that $\G_k (A)$ has an 
$\O\-$algebra structure and that $\R\G_k (A)$ is an ideal; then, for any $\chi\in \G_k (A)\,,$
let us denote by $\mu_\chi$ the multiplication by $\chi$ in $\G_k (A)\,.$
Moreover, denoting by  $\delta_a$ the characteristic $\K\-$valued function of $a\in A$  and by 
$A^*\i A$ the set of {\it generators\/} of $A\,,$ it is quite clear that $\{\delta_a\}_{a\in A^*}$
is a $\K\-$basis of $\K\otimes_\O\R\G_k (A)\,.$

\smallskip
In particular, the image $\bar C$ of $C$ in ${\rm Aut}(A)$ acts on $\K\otimes_\O\R\G_k (A)$ acting 
freely on this basis, and for any $\sigma\in C$ there is $\psi_\sigma\in {\rm Hom}(A,k^*)$ such that
$\sigma = \psi_\sigma\. \bar\sigma$ where $\bar\sigma$ denotes the image of $\sigma$ in ${\rm Aut}(A)$
(cf.~£3.2.15), so that for any $a\in A^*$ we get
$$\sigma (\delta_a) = (\mu_{\psi_\sigma}\circ \bar\sigma)(\delta_a) = \psi_\sigma\big(\bar\sigma (a)\big)\.\delta_{\bar\sigma (a)}
\eqno £3.2.17\phantom{.}$$
 where we identify the group $k^*$ with its canonical lifting to $\K^*\,;$ thus, if $O\i A^*$ is a $\bar C\-$orbit 
then the ideal $\I_O = \bigoplus_{a\in O} \K\.\delta_a$ of $\K\otimes_\O \G_k (A)$ is $C\-$stable
and, in order to prove equality~£3.2.13, it suffices to show that $(\I_O)^C$ has dimension one.
\eject

\smallskip
But, for any $\sigma,\tau\in C$ in ${\rm Aut}_{k^*}(\hat A)$ we get
$$\psi_{\tau\sigma}\.\overline{\tau\sigma} = \tau\sigma = (\psi_\tau\.\bar\tau)(\psi_\sigma\.\bar\sigma)
= \psi_\tau (\bar\tau\. \psi_\sigma\.\bar\tau^{-1})\.\overline{\tau\sigma}
\eqno £3.2.18\phantom{.}$$
which forces $\psi_{\tau \sigma} = \psi_\tau (\bar\tau\. \psi_\sigma\.\bar\tau^{-1})\,;$ then, choosing $a\in O\,,$
it is clear that the element
$$\chi = \sum_{\sigma\in C} \sigma(\delta_a) = \sum_{\sigma\in C}(\mu_{\psi_\sigma}\circ \bar\sigma)(\delta_a)
\eqno £3.2.19\phantom{.}$$
is invertible in $\I_O$ and, moreover, for any $\tau\in C\,,$ we have
$$\eqalign{\bar\tau (\chi) &= \sum_{\sigma\in C} \bar\tau\big((\mu_{\psi_\sigma}\circ \bar\sigma)(\delta_a)\big)
= \sum_{\sigma\in C} \big((\bar\tau\circ \mu_{\psi_\sigma}\circ\bar\tau^{-1})
\circ\overline{\tau\sigma}\big)(\delta_a)\cr
&= \sum_{\sigma\in C} (\mu_{\psi_\tau^{-1}}\circ \mu_{\psi_{\tau\sigma}}
\circ\overline{\tau\sigma})(\delta_a) = \psi_\tau^{-1}\chi\cr} 
\eqno £3.2.20;$$
consequently,  for any $a'\in O$ we get
$$\eqalign{(\mu_{\chi^{-1}}\circ\tau\circ \mu_\chi) (\delta_{a'}) &= (\mu_{\chi^{-1}}\circ\mu_{\psi_\tau})
\bar\tau(\chi\delta_{a'}) = \chi^{-1}\psi_\tau\bar\tau (\chi)\bar\tau (\delta_{a'})\cr
& =  \bar\tau (\delta_{a'})\cr} 
\eqno £3.2.21;$$
that is to say, the actions of $\tau$ and $\bar\tau$ on $\I_O$ are conjugate each other; since  
$(\I_O)^{\bar C}$ has cleraly dimension one, we are done.

\medskip
£3.3. We are ready to prove Proposition~£2.14, so  we assume that isomorphism~£2.14.1 holds; 
let $C$ be a cyclic subgroup of~${\rm Out}_{k^*} (\hat H)_c$ and $\theta$ an element in ${\rm Irr}_k (\hat L,c)\,;$
we already know that $C_\theta$ stabilizes a pair formed by $(R,g)\in \R$ and by a {\it projective\/} irreducible
character $\theta^*$ of $\skew4\hat{\tilde\F}_{\!(c,\hat L)}  (R)\,,$ and therefore it acts on $B_\theta$ (cf.~£2.11);
more precisely, $C_\theta$ acts on both $k^*\-$groups 
$\widehat B^{^\theta}_\theta$ and $\widehat B^{^{\theta^*}}_\theta$ and these actions lift its action 
on $B_\theta\,.$ Moreover, since $B_\theta$ is a cyclic $p'\-$group, we have $k^*\-$isomorphisms 
$$\widehat B^{^\theta}_\theta\cong B_\theta\times k^*\cong \widehat B^{^{\theta^*}}_\theta
\eqno £3.3.1\phantom{.}$$
lifting the identity on $B_\theta\,;$ setting $\hat B_\theta = B_\theta\times k^*\,,$ since ${\rm Aut}_{k^*}(\hat B)$ 
acts faithfully on $\G_k (\hat B)\,,$ if we assume that  the actions of $C_\theta$ on 
$\G_k (\hat H_\theta\!\mid\! \theta)$ and on $\G_k \big(\bar N_{\hat H_\theta}  (R,g)\!\mid\! \zeta_{\theta^*}\big)$ have the same kernel, then
it follows from isomorphisms £2.6.1 and~£2.11.2 that the images of~$C_\theta$ in 
${\rm Aut}_{k^*}(\widehat B^{^\theta}_\theta)$ and in ${\rm Aut}_{k^*}(\widehat B^{^{\theta^*}}_\theta)$
have the same order. Now, equality~2.14.2 follows from lemma~£3.2; we are done.

\bigskip
\noindent
{\bf £4. Proof of the first main result}
\bigskip
£4.1. In  [8,~Theorem~1.6],  for any block $b$ of any $k^*\-$group $\hat G$ of finite $k^*\-$quotient, 
 choosing a maximal Brauer $(b,\hat G)\-$pair $(P,e)$ and denoting by 
$\F_{\!(b,\hat G)}$ the  Frobenius $P\-$category of $(b,\hat G)$ (cf.~£1.5), we prove that the existence 
of an  $\O {\rm Out}_{k^*}(\hat G)_b\-$ module isomorphism
$$\G_k (\hat G,b)\cong \G_k (\F_{\!(b,\hat G)},\widehat\aut_{(\F_{\!(b,\hat G)})^{^{\rm sc}}})
 \eqno £4.1.1\phantom{.}$$
 \eject
 \noindent
  is equivalent to the existence of such an isomorphism for any block $c$ of any {\it almost simple\/} $k^*\-$group~$\hat H\,.$ Although not explicit, it is easy to check that, in all the steps of the proof,
if we only assume our isomorphisms defined over the corresponding $\K\-$extensions, we still obtain isomorphisms defined
over the $\K\-$extensions. Consequently, we may apply the following result to our present situation.

\bigskip
\noindent
{\bf Theorem~£4.2.} {\it Assume {\rm (SOSFG)}. If for any block 
$c$ of~positive defect of any almost simple $k^*\-$group $\hat H$ there is a $\K {\rm Out}_{k^*}(\hat H)_c\-$module isomorphism
$$\K\otimes_\O\G_k  (\hat H,c)\cong
\K\otimes_\O\G_k(\F_{\!(c,\hat H)},\widehat \aut_{(\F_{\!(c,\hat H)})^{^{\rm sc}}})
\eqno £4.2.1,$$
then for any block $b$ of any $k^*\-$group $\hat G$ of finite $k^*\-$quotient there is an  
$\K {\rm Out}_{k^*} (\hat G)_b\-$module isomorphism
$$\K\otimes_\O \G_k (\hat G,b)\cong  \K\otimes_\O \G_k(\F_{\!(b,\hat G)},
\widehat \aut_{(\F_{\!(b,\hat G)})^{^{\rm sc}}})
\eqno £4.2.2.$$\/}

\par
£4.3. With the same notation, it is clear that  isomorphism~£4.2.2 is equi-valent to the equality of the  ${\rm Out}_{k^*} (\hat G)_b\-$characters of both members which, as in~£2.1 above, amounts to saying that
  for any cyclic subgroup $C$ of ${\rm Out}_{k^*}(\hat G)_b$ we have
$${\rm rank}_\O \Big(\G_k(\F_{\!(b,\hat G)},
\widehat \aut_{(\F_{\!(b,\hat G)})^{^{\rm sc}}})^C\Big) =  {\rm rank}_\O\big(\G_k (\hat G,b)^C\big)
 \eqno £4.3.1\,.$$
 Consequently, in order to prove that Theorem~£4.2 implies Theorem~£1.11 it suffices to show that, 
 under the hypothesis in Theorem~£1.11, we have (cf.~£1.9)
 $$\eqalign{{\rm rank}_\O \Big(\G_k(\F_{\!(b,\hat G)},
&\,\widehat \aut_{(\F_{\!(b,\hat G)})^{^{\rm sc}}})^C\Big)\cr
& =   \sum_{(Q,f)} {\rm rank}_{\O} \Big(\G_k \big(\skew4\hat{\tilde\F}_{\!(b,\hat G)}  
(Q),o_{(Q,f)}\big)^{C_{(Q,f)}}\Big)\cr}
\eqno £4.3.2\phantom{.}$$
where $(Q,f)$ runs over a set of representatives contained in $(P,e)$ for the set of $C\-$orbits of 
$G\-$conjugacy classes of selfcentralizing  Brauer $(b,\hat G)\-$pairs  and, for such a $(Q,f)\,,$ we denote by 
$C_{(Q,f)}$ the stabilizer of the $G\-$conjugacy class of $(Q,f)$ in~$C\,.$ Indeed, our hypothesis in Theorem~£1.11 implies that any block $c$ of any almost-simple  $k^*\-$group  $\hat H$ fulfills equalities~£1.10.1; 
then, equalities~£4.3.2 show that the pair $(c,\hat H)$ also fulfills equalities~£4.3.1 and therefore isomorphism~£4.2.1 holds. At this point, Theorem~£4.2 implies that, for
any block $b$ of any $k^*\-$group $\hat G$ of finite $k^*\-$quotient, isomorphism~£4.2.2 holds
and therefore the pair $(b,\hat G)$ fulfills equalities~£4.3.1; finally, this time equalities~£4.3.2 show that
the pair $(b,\hat G)$ fulfills equalities~£1.10.1.

\medskip
£4.4. Note that, arguing by induction on $\vert G\vert\,,$ under the hypothesis of Theorem~£1.11
 we may assume that for any block $c$ of any $k^*\-$group $\hat H$ such that $\vert H\vert < \vert G\vert$
 and  any cyclic subgroup $D$ of ${\rm Out}_{k^*}(\hat H)_{c}$
 we have
 $$ {\rm rank}_\O\big(\G_k (\hat H,c)^D\big) = \sum_{(R,g)} {\rm rank}_{\O} \Big(\G_k 
 \big(\skew4\hat{\tilde\F}_{\!(c,\hat H)}  (R),o_{(R,g)}\big)^{D_{(R,g)}}\Big)
 \eqno £4.4.1\phantom{.}$$
 where $(R,g)$ runs over a set of representatives contained in a maximal Brauer $(c,\hat H)\-$pair 
 for the set of $D\-$orbits of  $H\-$conjugacy classes of selfcentralizing  Brauer $(c,\hat H)\-$pairs. But,
 as in £2.10 above, denoting by  $n_{(R,g)}$ the sum of all the blocks  of {\it defect zero\/} 
of~$k_*\bar N_{\hat H} (R,g)\bar g\,,$ it follows again from  [7,~Proposition~3.2 and~Theorem~3.7]
that we have a canonical isomorphism
$$ \G_k  \big(\bar N_{\hat H} (R,g),n_{(R,g)}\big) \cong \G_k  \big(\skew4\hat{\tilde\F}_{\!(c,\hat H)}  
(R),o_{(R,g)}\big)
\eqno £4.4.2;$$
actually, ${\rm Tr}_{\bar N_{\hat H} (R,g)}^{\bar N_{\hat H} (R)}(n_{(R,g)})$ is a  sum of blocks  
of {\it defect zero\/}  of~$\bar N_{\hat H} (R)$ and all the  blocks   of {\it defect zero\/}  of~$\bar N_{\hat H} (R)$
involved in ${\rm Br}_R(c)$ appear in these sums. Moreover, for any selfcentralizing  Brauer $(b,\hat G)\-$pair $(Q,f)$ we have $n_{(Q,f)}\not= 0$ only if 
$\Bbb O_p\big(\bar N_{\hat G} (Q,f)\big) = \{1\}\,;$ in this case,
 it is easily checked that any normal $p\-$subgroup $U$ of $\hat G$ is contained in $Q$ and, setting
 $\skew3\hat{\bar G} = \hat G/U$ and $\bar Q= Q/U$ we clearly  have $\bar N_{\skew3\hat{\bar G}}(\bar Q)
 \cong \bar N_{\hat G} (Q)\,.$ Consequently, if $\Bbb O_p(\hat G)\not = \{1\}$ then it is easily checked 
 from the induction hypothesis that the pair $(b,\hat G)$ also fulfills equality~£4.4.1

\medskip
£4.5. As a matter of fact, in the sequel it is more convenient to consider the {\it exterior\/} quotient
$\tilde\F_{\!(b,\hat G)}$ of $\F_{\!(b,\hat G)}$ [6,~1.3] formed by the same objects and by the morphisms
$\tilde\varphi\,\colon R\to Q$ where $\tilde\varphi$ denotes the $Q\-$conjugacy class of an  
$\F_{\!(b,\hat G)}\-$morphism $\varphi\,\colon R\to Q\,,$ the composition being induced by the composition in
$\F_{\!(b,\hat G)}\,;$ similarly, $(\tilde\F_{\!(b,\hat G)})^{^{\rm sc}}$ is the {\it full\/} subcategory
of $\tilde\F_{\!(b,\hat G)}$ determined by the set of selfcentralizing Brauer $(b,\hat G)\-$pairs contained
in $(P,e)\,.$ As in~£1.6 above, we consider  the  {\it proper category of 
$(\tilde\F_{\!(b,\hat G)})^{^{\rm sc}}\-$chains\/} $\ch^*\big((\tilde\F_{\!(b,\hat G)})^{^{\rm sc}}\big)$
and we still have the corresponding functor [6,~Proposition~A2.10]
$$\aut_{(\tilde\F_{\!(b,\hat G)})^{^{\rm sc}}} : \ch^*\big((\tilde\F_{\!(b,\hat G)})^{^{\rm sc}}\big)\too \Gr
\eqno £4.5.1.$$
In [6,~14.9] we show that this functor can also be lifted to an essentially unique functor
$$\widehat\aut_{(\tilde\F_{\!(b,\hat G)})^{^{\rm sc}}} : \ch^*\big((\tilde\F_{\!(b,\hat G)})^{^{\rm sc}}\big)
\too k^*\-\Gr
\eqno £4.5.2\phantom{.}$$
which composed with the canonical funtor
$$\ch^*\big((\F_{\!(b,\hat G)})^{^{\rm sc}}\big)\too \ch^*\big((\tilde\F_{\!(b,\hat G)})^{^{\rm sc}}\big)
\eqno £4.5.3\phantom{.}$$
\eject
\noindent
admits a natural map from $\widehat\aut_{(\F_{\!(b,\hat G)})^{^{\rm sc}}}$ (cf.~£1.7.1); then,
 for any  $(\tilde\F_{\!(b,\hat G)})^{^{\rm sc}}\-$ chain $\frak q\,\colon \Delta_n\to 
 (\tilde\F_{\!(b,\hat G)})^{^{\rm sc}}\,,$ we denote by $\skew4\hat{\tilde\F}_{\!(b,\hat G)}  (\frak q)$
 the image of $\frak q$ by the functor in~£4.5.2.

\medskip
£4.6. We actually will prove  equality~£4.3.2 in two steps; on the one hand, we will adapt our arguments in
the proof of [6, Corollary~14.32] in order to show that 
$$\eqalign{{\rm rank}_\O \Big(\G_k(\F_{\!(b,\hat G)},
  &\widehat \aut_{(\F_{\!(b,\hat G)})^{^{\rm sc}}})^C\Big)\cr
& = \sum_{\frak q\in \frak Q_{(b,\hat G)}} (-1)^{\ell(\frak q)} {\rm rank}_\O\Big(\G_k 
\big(\skew4\hat{\tilde \F}_{\!(b,\hat G)} (\frak q)\big)^{C_\frak q}\Big)\cr}
\eqno £4.6.1\phantom{.}$$
where $\frak Q_{(b,\hat G)}$ is a set of representatives, {\it fully normalized\/} in  $\F_{\!(b,\hat G)}$
(see £A7 below), for the  set of $\tilde\F_{\!(b,\hat G)}\-$isomorphism classes  of {\it regular\/} 
$(\tilde\F_{\!(b,\hat G)})^{^{\rm sc}}\-$chains (see £A6 below) and,
for such a~$\frak q\,\colon \Delta_n\to (\tilde\F_{\!(b,\hat G)})^{^{\rm sc}}\,,$ $C_\frak q$ denotes 
the stabilizer in $C$ of the isomorphism class of  $\frak q$ and we set $\ell (\frak q) = n\,.$
 On the other hand,  from Lemmas~A13 and~A14 below and our induction hypothesis we will prove that
$$\eqalign{\sum_{\frak q\in \frak Q_{(b,\hat G)}} (-1)^{\ell(\frak q)} 
{\rm rank}_\O &\Big(\G_k  \big(\skew4\hat{\tilde \F}_{\!(b,\hat G)} (\frak q)\big)^{C_\frak q}\Big)\cr
&= \sum_{(Q,f)} {\rm rank}_{\O} \Big(\G_k \big(\skew4\hat{\tilde\F}_{\!(b,\hat G)}  (Q),o_{(Q,f)}\big)^{C_{(Q,f)}}\Big)\cr}
\eqno £4.6.2\phantom{.}$$
where $(Q,f)$ runs over the same set of representatives as above (cf.~£1.9).

\medskip
£4.7. In the first step, we need some notation from [6,~Ch.~14]; for any $h\in \Bbb N - p\Bbb N\,,$ 
let us denote by $U_h$ the group of  $h\-$th roots of unity in $\O^*$ and by 
$({}^h \tilde\F_{\!(b,\hat G)})^{^{\rm sc}}$ the category where the objects are
the pairs $Q^\rho$ determined by a selfcentralizing Brauer $(b,\hat G)\-$pair  $(Q,f)$ contained in $(P,e)$ (cf.~£1.5) and 
by an injective group homomorphism $\rho\,\colon U_h\to \tilde\F_{\!(b,\hat G)}(Q)$ (cf.~£1.9.2), and where the morphisms  from another such a pair $R^\sigma$ to $Q^\rho$ are the $\tilde\F_{\!(b,\hat G)}\-$morphisms 
$\tilde\varphi\,\colon R\to Q$ such that,  for any $\xi\in U_h\,,$ we have
 $\sigma (\xi)\circ \tilde\varphi= \tilde\varphi\circ \rho (\xi)$ [6,~14.25]. Similarly, we denote by ${}^{U_h} \aleph$ the category of finite $U_h\-$sets --- namely, finite sets endowed with a $U_h\-$action --- and by [6,~14.21]
$$\F\!ct_{U_h} : {}^{U_h}\aleph\too \O\-\mod
\eqno £4.7.1\phantom{.}$$
the  {\it contravariant\/} functor  mapping any finite $U_h\-$set~$X$ on the $\O\-$module $\F\!ct_{U_h} (X,\O)$ of the
$\O\-$valued functions over~$X$ preserving the $U_h\-$actions --- $U_h$ acting on $\O$ by multiplication. Then, we consider the functor 
$$\frak s_h : ({}^h \tilde\F_{\!(b,\hat G)})^{^{\rm sc}}\too {}^{U_h} \aleph
\eqno £4.7.2\phantom{.}$$
provided by [6,~Proposition~14.28] and denote by ${}^\K \frak n_h$ the extension to $\K$ of the composed functor
$\F\!ct_{U_h}\circ \frak s_h\,.$
\eject

\medskip
£4.8. Now, it follows from [6,~Theorem~14.30] that, for any $n\ge 1\,,$ we have
$$\Bbb H^n \big( ({}^h \tilde\F_{\!(b,\hat G)})^{^{\rm sc}},{}^\K\frak n_h \big) = \{0\}
\eqno £4.8.1;$$
moreover, in [4,~A3.17] we consider the {\it stable cohomology\/} groups, denoted by 
$\Bbb H^n_* \big( ({}^h \tilde\F_{\!(b,\hat G)})^{^{\rm sc}},{}^\K\frak n_h \big)\,,$
computed from the {\it $n\-$cocycles\/} and the {\it $n\-$cobaundaries\/} which are ``stable''
by the obvious isomorphisms and, since we are working over the field $\K\,,$ in [6,~Propositions~A4.13]
we prove that, for any $n\in \Bbb N\,,$ we have
$$\Bbb H^n_* \big( ({}^h \tilde\F_{\!(b,\hat G)})^{^{\rm sc}},{}^\K\frak n_h \big)\cong
\Bbb H^n \big( ({}^h \tilde\F_{\!(b,\hat G)})^{^{\rm sc}},{}^\K\frak n_h \big)
\eqno £4.8.2.$$
Finally, it is quite clear that the category $({}^h \tilde\F_{\!(b,\hat G)})^{^{\rm sc}}$ fulfills
the condition [6,~A5.1.1] and we can consider the  {\it regular\/} 
$ ({}^h \tilde\F_{\!(b,\hat G)})^{^{\rm sc}}\-$chains, namely 
the $ ({}^h \tilde\F_{\!(b,\hat G)})^{^{\rm sc}}\-$chains $\frak q^\eta\,\colon \Delta_n\too 
({}^h \tilde\F_{\!(b,\hat G)})^{^{\rm sc}}$ [6,~Proposition~14.27] such that $\frak q^\eta(i\!-\!1\bullet i)$ 
is {\it not\/} an isomorphism for any $1\le i\le n$ [6,~A5.2]; then, in [6,~Proposition~A4.7] we show that the groups $\Bbb H^n_* \big( 
({}^h \tilde\F_{\!(b,\hat G)})^{^{\rm sc}},{}^\K\frak n_h \big)$ can be computed from
the {\it regular\/} $ ({}^h \tilde\F_{\!(b,\hat G)})^{^{\rm sc}}\-$chains, namely that,  for any $n\in \Bbb N\,,$ we have
$$\Bbb H^n_* \big( ({}^h \tilde\F_{\!(b,\hat G)})^{^{\rm sc}},{}^\K\frak n_h \big)\cong
\Bbb H^n_{\rm r} \big( ({}^h \tilde\F_{\!(b,\hat G)})^{^{\rm sc}},{}^\K\frak n_h \big)
\eqno £4.8.3.$$

\medskip
£4.9. In conclusion,  for any $n\ge 1\,,$ we have
$$\Bbb H^n_{\rm r} \big( ({}^h \tilde\F_{\!(b,\hat G)})^{^{\rm sc}},{}^\K\frak n_h \big) = \{0\}
\eqno £4.9.1;$$
that is to say, for any $n\in \Bbb N$ setting 
$$\Bbb C^n_{\rm r} = \prod_{\frak q^\eta} \K\otimes_\O\F\!ct_{U_h} 
\Big(\frak s_h \big(\frak q^\eta (0)\big),\O\Big)^{\tilde\F_{\!(b,\hat G)}(\frak q)_\eta}
\eqno £4.9.2\phantom{.}$$
where $\frak q$ runs over a set of representatives for the set of isomorphism
classes of {\it regular\/} $ (\tilde\F_{\!(b,\hat G)})^{^{\rm sc}}\-$chains, $\eta\,\colon U_h\to \tilde\F (\frak q)$
runs over the set of injective group homomorphisms [6,~Proposition~14.27] and $\tilde\F_{\!(b,\hat G)}(\frak q)_\eta$
denotes the stabilizer of $\eta$ in $\tilde\F_{\!(b,\hat G)}(\frak q)\,,$ we have a {\it finite\/} exact sequence
$$0\too \Bbb H^0 \big( ({}^h \tilde\F_{\!(b,\hat G)})^{^{\rm sc}},{}^\K\frak n_h\big)
\too  \Bbb C^0_{\rm r}\too \dots\too \Bbb C^{n}_{\rm r}\too \dots
\eqno £4.9.3.$$
But, since we are working over $\K\,,$   for any cyclic subgroup $C$ of ${\rm Out}_{k^*}(\hat G)_b$ we still have the {\it finite\/} exact sequence of $C\-$fixed points
$$0\too \Bbb H^0 \big( ({}^h \tilde\F_{\!(b,\hat G)})^{^{\rm sc}},{}^\K\frak n_h\big)^C
\too  (\Bbb C^0_{\rm r})^C\too\dots\too (\Bbb C^{n}_{\rm r})^C\too \dots
\eqno £4.9.4.$$
Consequently, we still get
$$\eqalign{&{\rm dim}_\K \Big(\Bbb H^0 \big( ({}^h \tilde\F_{\!(b,\hat G)})^{^{\rm sc}},{}^\K\frak n_h\big)^C\Big)\cr
& = \sum_{\frak q^\eta} (-1)^{\ell(\frak q^\eta)}\,{\rm dim\,}_\K 
\bigg(\K\otimes_\O\F\!ct_{U_h} \Big(\frak s_h \big(\frak q^\eta (0)\big),\O\Big)^{\tilde\F_{\!(b,\hat G)}(\frak q)_\eta
\rtimes C_{\tilde\frak q^\eta}}\bigg)\cr}
\eqno £4.9.5\phantom{.}$$
where $\frak q^\eta$ runs over a set of representatives for the isomorphism
classes  of {\it regular\/} $ ({}^h \tilde\F_{\!(b,\hat G)})^{^{\rm sc}}\-$chains  [6,~Proposition~14.27] and 
$C_{\frak q^\eta}$ denotes the stabilizer in $C$ of the isomorphism class of $\frak q^\eta\,.$

\medskip
£4.10. Finally, on the one hand, it follows from [6,~14.28.3] that we have
$$\eqalign{{\rm rank}_\O \Big(\G_k(\F_{\!(b,\hat G)},&\,\widehat \aut_{(\F_{\!(b,\hat G)})^{^{\rm sc}}})^C\Big)\cr
&= \sum_{h\in \Bbb N - p\Bbb N} {\rm dim}_\K \Big(\Bbb H^0 
\big(({}^h \tilde\F_{\!(b,\hat G)})^{^{\rm sc}},{}^\K\frak n_h\big)^C\Big)\cr}
\eqno £4.10.1.$$
 On the other hand, for any {\it regular\/} $ (\tilde\F_{\!(b,\hat G)})^{^{\rm sc}}\-$chain $\frak q\,,$
 setting $\hat F = \skew4\hat{\tilde \F}_{\!(b,\hat G)} (\frak q)$ and $F = \tilde \F_{\!(b,\hat G)} (\frak q)\,,$ 
it follows from [6,~14.15.3] that we have
 $$\G_k (\hat F)^{C_\frak q}\cong \bigg(\bigoplus_{h\in \Bbb N -p\Bbb N}\,
 \bigoplus_{\eta\in {\rm Mon}(U_h,F)} \F\!ct_{U_h} \Big((\varpi_{h,\hat F})^{-1}(\eta)  ,
\O\Big)\bigg)^{F\rtimes C_{\frak q}} 
\eqno £4.10.2\phantom{.}$$
where $C_\frak q$ denotes the stabilizer in $C$  of the isomorphism class of $\frak q$ and, for any 
$h\in \Bbb N -p\Bbb N\,,$  setting $\hat U_h = U_h\times k^*$ we respectively denote by 
$${\rm Mon}\big(U_h,\tilde\F_{\!(b,\hat G)}(\frak q)\big)\qq {\rm Mon}_{k^*}
\big(\hat U_h,\skew4\hat{\tilde\F}_{\!(b,\hat G)}(\frak q)\big)
\eqno £4.10.3\phantom{.}$$
 the sets of injective group and $k^*\-$group homomorphisms from $U_h$ to $\tilde\F_{\!(b,\hat G)}(\frak q)$ and
 from $\hat U_h$ to $\skew4\hat{\tilde\F}_{\!(b,\hat G)}(\frak q)\,,$ and by
$$\varpi_{h,\skew4\hat{\tilde \F}_{\!(b,\hat G)} (\frak q)} : {\rm Mon}_{k^*}
\big(\hat U_h,\skew4\hat{\tilde\F}_{\!(b,\hat G)}(\frak q)\big)\too
{\rm Mon}\big(U_h,\tilde\F_{\!(b,\hat G)}(\frak q)\big)
\eqno £4.10.4\phantom{.}$$
the canonical map.

\medskip
£4.11. But, by the very definition of the functor $\frak s_h$ in [6,~Proposition~14.28], for any {\it regular\/} 
$ ({}^h \tilde\F_{\!(b,\hat G)})^{^{\rm sc}}\-$chain $\frak q^\eta$ we have
$$\frak s_h \big(\frak q^\eta (0)\big) = (\varpi_{h,\skew4\hat{\tilde\F}_{\!(b,\hat G)}
(\frak q (0))})^{-1}(\iota_0^\frak q\circ\eta)\cong (\varpi_{h,\skew4\hat{\tilde\F}_{\!(b,\hat G)}
(\frak q )})^{-1}(\eta)
\eqno  £4.11.1\phantom{.}$$
where $\iota_0^\frak q\,\colon \tilde\F (\frak q)\to \tilde\F\big(\frak q (0)\big)$ is the structural map
[6,~14.26]. Consequently, isomorphism~£4.10.2 becomes
$$\eqalign{&\G_k \big(\hat \F_{\!(b,\hat G)} (\frak q)\big)^{C_\frak q}\cr
&\cong \bigg(\bigoplus_{h\in \Bbb N -p\Bbb N}\,\bigoplus_{\eta\in {\rm Mon}(U_h,\tilde\F_{\!(b,\hat G)}
(\frak q))} \F\!ct_{U_h} \Big(\frak s_k \big(\frak q^\eta (0)\big),
\O\Big)\bigg)^{\tilde\F_{\!(b,\hat G)}(\frak q )\rtimes C_{\frak q}}\cr
&\cong \bigoplus_{h\in \Bbb N -p\Bbb N}\,\bigoplus_{\eta\in {\rm Mon}(U_h,\tilde\F_{\!(b,\hat G)}(\frak q))}
 \F\!ct_{U_h} \Big(\frak s_k \big(\frak q^\eta (0)\big),\O\Big)^{\tilde\F_{\!(b,\hat G)}(\frak q)_\eta
\rtimes C_{\frak q^\eta}}\cr}
\eqno £4.11.2\phantom{.}$$
In conclusion, the sum of all the equalities~£4.9.5 when $h$ runs over $\Bbb N -p\Bbb N$ yields
 $$\eqalign{{\rm rank}_\O \Big(\G_k(\F_{\!(b,\hat G)},
 &\widehat \aut_{(\F_{\!(b,\hat G)})^{^{\rm sc}}})^C\Big)\cr
& = \sum_{\frak q\in \frak Q_{(b,\hat G)}} (-1)^{\ell(\frak q)} {\rm rank}_\O\Big(\G_k 
\big(\skew4\hat{\tilde\F}_{\!(b,\hat G)} (\frak q)\big)^{C_\frak q}\Big)\cr}
\eqno £4.11.3,$$
proving equality~£4.6.1.

\medskip
£4.12. Before going further, note that in the right-hand member many pairs of terms in this sum cancel each other;
more precisely, it easily follows from Lemma~A13 below that, in this sum, we can replace the set $\frak Q_{(b,\hat G)}$
by the image in $\ch^*\big((\tilde\F_{\!(b,\hat G)})^{^{\rm sc}}\big)$ of  a set of representatives for
the isomorphism classes in the intersection  $\frak N^{\F_{\!(b,\hat G)}}$ and, similarly, from Lemma~A14 below that we still can replace this set by  the image in  $\ch^*\big((\tilde\F_{\!(b,\hat G)})^{^{\rm sc}}\big)$ of  a set of representatives for the isomorphism classes in the  intersection~$\frak R^{\F_{\!(b,\hat G)}}\,;$ explicitly, we get
 $$\eqalign{{\rm rank}_\O \Big(\G_k(\F_{\!(b,\hat G)},
 &\widehat \aut_{(\F_{\!(b,\hat G)})^{^{\rm sc}}})^C\Big)\cr
& = \sum_{\check\frak q\in \frak R^{\F_{\!(b,\hat G)}}} (-1)^{\ell(\frak q)} {\rm rank}_\O\Big(\G_k 
\big(\skew4\hat{\tilde\F}_{\!(b,\hat G)} (\frak q)\big)^{C_{\frak q}}\Big)\cr}
\eqno £4.12.1\phantom{.}$$
where, for any $(\tilde\F_{\!(b,\hat G)})^{^{\rm sc}}\-$chain $\frak q\,,$ we denote by $\check\frak q$ the isomorphism
class of any lifting of $\frak q$ to an $(\F_{\!(b,\hat G)})^{^{\rm sc}}\-$chain.

\medskip
£4.13. On the other hand, let $\frak q\,\colon \Delta_n\to (\tilde\F_{\!(b,\hat G)})^{^{\rm sc}}$ be a {\it regular\/} 
$(\tilde\F_{\!(b,\hat G)})^{^{\rm sc}}\!\-$ chain {\it fully normalized\/} in $\tilde\F_{\!(b,\hat G)}$ 
(see~A6 and~A7 below where we replace $\F_{\!(b,\hat G)}$ by~$\tilde\F_{\!(b,\hat G)}\,$)  and consider the corresponding normalizer $N_{\hat G} (\frak q)\,;$ since the Brauer $(b,\hat G)\-$pair  $(Q_\frak q,b_\frak q)$ determined by $\frak q(n) = Q_\frak q$ and by the condition  $(Q_\frak q,b_\frak q)\i (P,e)$ is selfcentralizing (cf.~£1.5), $b_\frak q $ is
actually   a {\it nilpotent\/} block of $C_{\hat G}(Q_\frak q)\,;$ moreover, it is well-known and easily checked that
$b_\frak q$ remains a block of~$N_{\hat G} (\frak q)\,.$ Then, it follows from [10,~Theorem~3.5 and~Corollary~3.15]  (see also [2,~Theorems~1.8 and 1.12]) that there exists a suitable $k^*\-$group $\hat L_{\F_{\!(b,\hat G)}} (\frak q)\,,$
containing $N_P (\frak q)$ and admitting the exact sequence
$$1\too \frak q(0)\too \hat L_{\F_{\!(b,\hat G)}} (\frak q)\too \skew4\hat{\tilde\F}_{\!(b,\hat G)} (\frak q)\too 1
\eqno £4.13.1,$$
such that we have canonical isomorphisms
$$\G_k (N_{\hat G} (\frak q),b_\frak q)\cong \G_k \big(\hat L_{\F_{\!(b,\hat G)}} (\frak q)\big)
\cong \G_k\big(\skew4\hat{\tilde\F}_{\!(b,\hat G)} (\frak q)\big)
\eqno £4.13.2.$$

\medskip
£4.14.  Actually, it is easily checked from the corresponding definitions that the $k^*\-$quotient 
 $L_{\F_{\!(b,\hat G)}} (\frak q)$ of  $\hat L_{\F_{\!(b,\hat G)}} (\frak q)$ coincides with the  {\it $\F_{\!(b,\hat G)}\-$localizer\/} $L_{\F_{\!(b,\hat G)}} (\frak q)$ of $\frak q$
[6,~Theorem~18.6]; in particular, from [6,~Remark~18.7] it is easy to check that we have the equivalence of categories
$$N_{\F_{\!(b,\hat G)}} (\frak q) \cong \F_{\! L_{\F_{\!(b,\hat G)}} (\frak q)}
\eqno £4.14.1\phantom{.}$$
where $\F_{\! L_{\F_{\!(b,\hat G)}} (\frak q)}$ denotes the {\it Frobenius category\/} associated with the group
$ L_{\F_{\!(b,\hat G)}} (\frak q)$ [6,~1.8]; note that $1$ is the unique block of $\hat L_{\F_{\!(b,\hat G)}} (\frak q)$
and that we have
$$\F_{\! L_{\F_{\!(b,\hat G)}} (\frak q)} = \F_{\!(1,\hat L_{\F_{\!(b,\hat G)}} (\frak q))}
\eqno £4.14.2.$$
Moreover, it follows from [6,~Corollary~3.6] that we also have an equivalence of categories
$$N_{\F_{\!(b,\hat G)}} (\frak q) \cong \F_{\!(b_\frak q, N_{\hat G}(\frak q))}
\eqno £4.14.3,$$
so that  the blocks $b_\frak q$ of $N_{\hat G}(\frak q)$
and $1$ of $\hat L_{\F_{\!(b,\hat G)}} (\frak q)$ have the same {\it Frobenius $N_P (\frak q)\-$category\/}.

\medskip
£4.15. More precisely, it follows from [10,~Corollary~3.15] that the blocks $b_\frak q$ of $N_{\hat G}(\frak q)$
and $1$ of $\hat L_{\F_{\!(b,\hat G)}} (\frak q)$ are {\it basic Morita equivalent\/} [5,~7.1]; in par-ticular,
it follows from [5,~7.6] that the corresponding functors 
$$\widehat\aut_{ (\F_{\!(b_\frak q, N_{\hat G}(\frak q))})^{^{\rm sc}}}\qq 
\widehat\aut_{(\F_{\!(1,\hat L_{\F_{\!(b,\hat G)}} (\frak q))})^{^{\rm sc}}}
\eqno £4.15.1\phantom{.}$$
are isomorphic; hence, denoting by ${\rm Out}_{k^*}(\hat G)_{b, \frak q}$ the stabilizer of the $G\-$conju-gacy 
class of $\frak q$ in ${\rm Out}_{k^*}(\hat G)_{b}\,,$ we also have an 
$\O {\rm Out}_{k^*}(\hat G)_{b, \frak q}\-$isomorphism
$$\eqalign{\G_k (\F_{\!(b_\frak q, N_{\hat G}(\frak q))},&\widehat\aut_{ (\F_{\!(b_\frak q, N_{\hat G}
(\frak q))})^{^{\rm sc}}})\cr
&\cong \G_k \big(\F_{\!(1,\hat L_{\F_{\!(b,\hat G)}} (\frak q))}, 
\widehat\aut_{(\F_{\!(1,\hat L_{\F_{\!(b,\hat G)}} (\frak q))})^{^{\rm sc}}}\big)\cr}
\eqno £4.15.2.$$
But, by the very definition of {\it inverse limit\/}, it is quite clear that we have an 
$\O {\rm Out}_{k^*}\big(\hat L_{\F_{\!(b,\hat G)}} (\frak q)\big)\-$isomorphism (cf.~£1.7.2)
$$\G_k \big(\F_{\!(1,\hat L_{\F_{\!(b,\hat G)}} (\frak q))}, 
\widehat\aut_{(\F_{\!(1,\hat L_{\F_{\!(b,\hat G)}} (\frak q))})^{^{\rm sc}}}\big)\cong 
\G_k \big(\hat L_{\F_{\!(b,\hat G)}} (\frak q)\big)
\eqno £4.15.3.$$
Consequently, from isomorphisms~£4.13.2 we obtain    $\O {\rm Out}_{k^*}(\hat G)_{b, \frak q}\-$isomor-phisms
$$\G_k (\F_{\!(b_\frak q, N_{\hat G}(\frak q))},\widehat\aut_{ (\F_{\!(b_\frak q, N_{\hat G}(\frak q))})^{^{\rm sc}}})\cong \G_k (N_{\hat G} (\frak q),b_\frak q)\cong \G_k\big(\skew4\hat{\tilde\F}_{\!(b,\hat G)} (\frak q)\big)
\eqno £4.15.4.$$

\medskip
£4.16.  In particular, if we have $\hat G = N_{\hat G}(\frak q)\,,$ either $\frak q (0) = \{1\}$ which forces
$P = \{1\}\,,$ so that $b$ is a block of {\it defect zero\/} and equality~£4.6.2 is {\it tautologically\/} true, 
or we have $\Bbb O_p (\hat G)\not=  \{1\}$ and therefore the pair $(b,\hat G)$ fulfills equality~£4.4.2 (cf.~£4.4);
in this case, we get equality~£4.6.2 from the left-hand isomorphism in~£4.15.4. Otherwise, we have
$\hat G \not= N_{\hat G}(\frak q)$ for any  {\it regular\/}  $(\tilde\F_{\!(b,\hat G)})^{^{\rm sc}}\!\-$ chain
$\frak q$  {\it fully normalized\/} in $\tilde\F_{\!(b,\hat G)}$ and therefore it follows from our induction hypothesis (cf.~£4.4) that  for any cyclic subgroup $C$ of 
${\rm Out}_{k^*}(\hat G)_b\,,$ denoting by $C_\frak q$ the stabilizer in $C$ of the isomorphism class of $\frak q\,,$
 we get
$$\eqalign{{\rm rank}_\O\Big(\G_k&\big(\skew4\hat{\tilde\F}_{\!(b,\hat G)} (\frak q)\big)^{C_\frak q}\Big)\cr 
 &=  \sum_{(Q,f)\in \Q_\frak q}{\rm rang}_\O \Big(  \G_k \big(\skew4\hat{\tilde\F}_{\!(b_\frak q,N_{\hat G}
 (\frak q))}  (Q),o_{(Q,f)}\big)^{C_{\frak q,(Q,f)}}\Big)\cr}
\eqno £4.16.1.$$
where $\Q_\frak q$ denotes a set of representatives, contained in a maximal Brauer 
$(b_\frak q,N_{\hat G} (\frak q))\-$pair, for the set of  $N_G (\frak q)\-$conjugacy classes of selfcentralizing Brauer $(b_\frak q,N_{\hat G} (\frak q))\-$pairs
and, for any $(Q,f)\in \Q_\frak q\,,$ $C_{\frak q,(Q,f)}$ denotes the stabilizer in $C_\frak q$ of the 
$N_G (\frak q)\-$conjugacy class of $(Q,f)\,;$ note that we have $o_{(Q,f)}\not= 0$ only if 
$$\Bbb O_p\big(\skew4\hat{\tilde\F}_{\!(b_\frak q,N_{\hat G} (\frak q))}  (Q)\big) = \{1\}
\eqno £4.16.2,$$
so that only if $Q$ is {\it $\F_{\!(b_\frak q,N_{\hat G} (\frak q))}\-$radical\/} (see £A4 below).

\medskip
£4.17 At this point, in order to prove equality~£4.6.2, we have to compute the double sum  (cf.~£4.12.1)
$$\sum_{\check\frak q\in \frak R^{\F_{\!(b,\hat G)}}}
\sum_{(R,g)\in \R_\frak q}  (-1)^{\ell(\frak q)}  {\rm rang}_\O 
\Big(  \G_k \big(\skew4\hat{\tilde\F}_{\!(b_\frak q,N_{\hat G} (\frak q))}  (R),o_{(R,g)}\big)^{C_{\frak q,(R,g)}}\Big)
\eqno £4.17.1.$$
where, for any $\check\frak q\in \frak R^{\F_{\!(b,\hat G)}}\,,$ $\R_\frak q$ is the subset of $R\in \Q_\frak q$
which are {\it $\F_{\!(b_\frak q,N_{\hat G} (\frak q))}\-$ radical\/}. But, for such a $\frak q\,\colon \Delta_n\to  (\F_{\!(b,\hat G)})^{^{\rm sc}}\,,$ $\frak  q(n)$ is actually a normal~$p\-$sub-group of
 $N_{\hat G} (\frak q)$ and thus, according to Lemma~£A5 below, $R$ contains~$\frak q (n)\,;$ thus, if $R\cong \frak q(0)$ then $n = 0\,;$ otherwise, either $R\not= \frak q (n)$ and we consider\break
\eject
\noindent 
 the
{\it regular\/} $(\F_{\!(b,\hat G)})^{^{\rm sc}}\-$chain $\frak q^\varpi\,\colon \Delta_{n+1}\to 
(\F_{\!(b,\hat G)})^{^{\rm sc}}$ extending $\frak q$ and mapping $n+1$ on $R$ and the 
$\Delta_n\-$morphism $n\bullet n+1$  on the corresponding 
inclusion map, or $R= \frak q (n)$  and we can consider the  restriction $\frak q^\varpi\,\colon \Delta_{n-1}
\to (\F_{\!(b,\hat G)})^{^{\rm sc}}$ of $\frak q\,.$ In both cases, $R$ remains an
{\it $N_{\F_{\!(b,\hat G)}} (\frak q^\varpi)\-$radical\/} subgroup of $N_P (\frak q^\varpi)$ and we clearly have
$$\eqalign{\Big(\G_k \big(\skew4\hat{\tilde\F}_{\!(b_{\frak q^\varpi},N_{\hat G} (\frak q^\varpi))}  
&(R),o_{(R,g} \big)\Big)^{C_{\frak q^\varpi,(R,g)}}\cr
& \cong \Big(\G_k \big(\skew4\hat{\tilde\F}_{\!(b_\frak q,N_{\hat G} (\frak q))}  (R),
o_{(R,g)}\big)\Big)^{C_{\frak q,(R,g)}}\cr}
\eqno £4.17.2;$$
moreover, we obviously have $(\frak q^\varpi)^\varpi = \frak q\,.$ In conclusion, since we clearly have 
$\vert \ell (\frak q^\varpi)- \ell (\frak q)\vert = 1\,,$ the double sum~£4.17.1 becomes 
$$\sum_{(R,g)} {\rm rank}_\O \Big(\G_k \big(\skew4\hat{\tilde\F}_{\!(b,\hat G)}  (R),o_{(R,g)}\big)^{C_{(R,g)}}\Big)
\eqno £4.17.3\phantom{.}$$
where $(R,g)$ runs over a set of representatives contained in $(P,e)$ for the set of 
$G\-$conjugacy classes of selfcentralizing Brauer $(b,\hat G)\-$pairs such that $R$ is a fully normalized 
{\it $\F_{\!(b,\hat G)}\-$radical\/} subgroup of $P\,,$ proving equality~£4.6.2. We are done.

\bigskip
\noindent
{\bf Appendix: Radical functions  over folded Frobenius categories}
\bigskip

£A1. The contents of this Appendix rises from [12] where Jacques Th\'eve-naz adopts the old point of view
consisting on that, in a finite group $G\,,$ the word ``local'' is synonymous  of ``concerning the family of normalizers of nontrivial  $p\-$subgroups''. Here we exhibit what seems a more adequate framework, involving Frobenius categories.
As a matter of fact, our arguments are useful in the proof of our main result above.
Let $P$ be a finite $p\-$group and denote  by  $\frak i\Gr$ the category formed by the finite groups and by the injective group  homomorphisms, and  by 
$\F_{\!P}$ the subcategory of  $\frak i\Gr$ where the objects are all the  subgroups
 of $P$ and the morphisms are the group homomorphisms induced by conjugation by elements of $P\,.$

 \medskip
£A2. Recall that a {\it Frobenius  $P\-$category\/} $\F$ is a subcategory 
of $\frak i\Gr$ containing $\F_{\!P}$ where the objects are all the  subgroups of $P$
and the morphisms fulfill the following three conditions [6, 2.8 and Proposition~2.11]
\smallskip
\noindent
£A2.1\quad {\it For any subgroup $Q$ of $P$ the inclusion functor $(\F)_Q\to 
(\frak i\Gr)_Q$ is full.\/}
\smallskip
\noindent
£A2.2\quad {\it $\F_P (P)$ is a Sylow $p\-$subgroup of $\F (P)\,.$\/}
\smallskip
\noindent
£A2.3\quad {\it If $Q$ is a subgroup of $P$ fulfilling $\xi \big(C_P (Q)\big)
= C_P\big(\xi (Q)\big)$ for any $\F\-$morphism $\xi\,\colon Q\.C_P(Q)\to P\,,$ if $\varphi\,\colon Q\to P$
is an $\F\-$morphism  and if $R$ is a subgroup of $N_P\big(\varphi(Q)\big)$ containing $\varphi (Q)$ such that $\F_P(Q)$ contains the action of $\F_R \big(\varphi(Q)\big)$ over $Q$ via $\varphi\,,$ then there is an $\F\-$morphism 
$\zeta\,\colon R\to P$ fulfilling $\zeta\big(\varphi (u)\big) = u$ for any $u\in Q\,.$\/}
\eject
\smallskip
\noindent
As in [6,~1.2], for any pair of subgroups $Q$ and $R$ of $P\,,$ we denote by $\F (Q,R)$ the set of $\F\-$morphisms from $Q$ to $R$ and set $\F (Q) = \F (Q,Q)\,;$ moreover, recall
that, for any category $\frak C$ and any $\frak C\-$object $C\,,$ $\frak C_C$ 
(or $(\frak C)_C$ to avoid confusion) denotes the category of ``$\frak C\-$morphisms
to $C$'' [6, 1.7].

\medskip
£A3. Given a Frobenius $P\-$category $\F\,,$ a subgroup $Q$ of $P$ and
a subgroup $K$ of the group ${\rm Aut}(Q)$ of automorphisms of $Q\,,$
we say that $Q$ is {\it fully $K\-$normalized\/} in $\F$ if we have [6, 2.6]
$$\xi \big(N_P^K (Q)\big) = N_P^{\,{}^\xi \!K}\big(\xi (Q)\big)
\eqno £A3.1$$
for any $\F\-$morphism $\xi\,\colon Q\.N_P^K (Q)\to P\,,$ where $N_P^K (Q)$ is the converse image of $K$ in $N_P (Q)$ {\it via\/} the canonical group homomorphism
$N_P (Q)\to {\rm Aut}(Q)$ and ${}^\xi \! K$ is the image of $K$ in 
${\rm Aut}\big(\xi (Q)\big)$ {\it via\/} $\xi\,.$ Recall that if $Q$ is fully 
$K\-$normalized in $\F$ then we have a new Frobenius $N_P^K(Q)\-$category 
$N_\F^K(Q)$ where, for any pair of subgroups $R$ and $T$ of $N_P^K (Q)\,,$
$\big(N_\F^K(Q)\big)(R,T)$ is the set of group homomorphisms {\it from $T$ to $R$\/} induced by the $\F\-$morphisms $\psi\,\colon Q\.T\to Q\.R$ which stabilize $Q$ and induce on it an element of $K$ [6, 2.14 and Proposition~2.16].

\medskip
£A4. We say that a subgroup $Q$ of $P$ is {\it $\F\-$selfcentralizing\/} if we have
$$C_P\big(\varphi (Q))\i \varphi (Q)
\eqno £A4.1\phantom{.}$$
 for any $\varphi \in \F (P,Q)\,,$ and we denote by $\F^{^{\rm sc}}$ the full subcategory of $\F$ over the set of $\F\-$selfcentralizing subgroups of $P\,.$ From the case of the Frobenius $P\-$categories {\it associated with a block of a finite group\/}, we know that it only makes sense to consider central $k^*\-$extensions of $\F (Q)$ whenever $Q$ is 
$\F\-$selfcentralizing [6,~7.4]; but, if $U$ is a subgroup of $P$ fully $K\-$normalized in~$\F$
for some subgroup $K$ of ${\rm Aut}(U)\,,$ a $N_{\F}^K (U)\-$selfcentralizing subgroup
of $N_P (Q)$ need not be  $\F\-$selfcentralizing, which is a handicap when comparing
choices of central $k^*\-$extensions in $\F$ and in $N_F^K (U)\,.$
 In order to overcome this difficulty, we consider the {\it $\F\-$radical\/}
subgroups of~$P\,;$ we say that a subgroup $R$ of $P$ is {\it $\F\-$radical\/} if it is
$\F\-$selfcentralizing and we have
$${\bf O}_p\big(\tilde\F (R)\big) = \{1\}
\eqno £A4.2\phantom{.}$$
where   $\tilde\F (R) = \F (R)/\F_R (R)$ [6, 1.3]; we denote by $\F^{^{\rm rd}}$ the full subcategory of $\F$ over the set of $\F\-$radical subgroups of $P\,.$

\bigskip
\noindent
{\bf Lemma~£A5}\phantom{.} {\it Let $\F$ be a Frobenius $P\-$category, $U$
a subgroup of $P$ and $K$ a subgroup of ${\rm Aut}(U)$ containing ${\rm Int}(U)\,.$
If $U$ is fully $K\-$normalized in $\F$ then any $N_\F^K(U)\-$radical subgroup $R$ of $N_P^K(U)$ contains $U$ and, in particular, it is $\F\-$selfcentralizing.\/}

\medskip
\noindent
{\bf Proof:} It is quite clear that the image of $N_{U\.R}(R)$ in $\big(N_\F^K(U)\big)(R)$
is a normal $p\-$subgroup and therefore it is contained in ${\bf O}_p\Big(\big(N_\F^K(U)\big)(R)\Big)\,,$ so
that $N_{U\.R}(R) = R$ which forces $U\.R = R\,.$
\eject 

\smallskip
Moreover,
for any $\F\-$morphism $\psi\,\colon R\to P\,,$ it is clear that $\psi (U)$ is a normal 
subgroup of $\psi (R)\.C_P\big(\psi (R)\big)$ and therefore, since $U$ is also fully centralized in $\F$ [6, Proposition~2.12], it follows from £A2.3 that there is an 
$\F\-$morphism 
$$\zeta  : \psi (R)\.C_P\big(\psi (R)\big)\too P
\eqno £A5.1\phantom{.}$$
fulfilling $\zeta \big(\psi (u)\big) = u$ for any $u\in U\,,$ so that the group homomorphism from $R$ to $N_P^K(U)$ mapping $v\in R$ on $\zeta\big(\psi (v)\big)$
is a $N_\F^K(U)\-$morphism; in particular, $\zeta\big(\psi (R)\big)$ is also 
$N_\F^K(U)\-$selfcentralizing and therefore we get
$$\zeta \Big(C_P\big(\psi (R)\big)\Big)\i \zeta\big(\psi (R)\big)
\eqno £A5.2\phantom{.}$$
which forces $C_P\big(\psi (R)\big)\i \psi (R)\,.$ We are done.

\medskip
£A6. Here, we have to deal with  {\it $\F^{^{\rm sc}}\-$chains\/} and {\it coherent\/} choices of central $k^*\-$extensions for the $\F^{^{\rm sc}}\!\-$automorphism groups. Recall that we call {\it $\F^{^{\rm sc}}\-$chain\/} any functor 
$\frak q\,\colon \Delta_n\to \F^{^{\rm sc}}$ where the $n\-$simplex $\Delta_n$ is considered as a category with the morphisms defined by the order relation [6,~A2.2];
let us call $n$ the {\it length\/} of $\frak q$ and
 set $n = \ell(\frak q)\,;$ recall that $\frak q$ is {\it regular\/} if 
 $\frak q(i\!-\!1\bullet i)$ is {\it not\/} an isomorphism
 for any $1\le i\le n$ [6,~A5.2]. Then, we consider the category 
$\ch^*(\F^{^{\rm sc}})$ where the objects are all the $\F^{^{\rm sc}}\-$chains 
$\frak q$ and the morphisms from $\frak q\,\colon \Delta_n\to 
\F^{^{\rm sc}}$ to another $\F^{^{\rm sc}}\-$chain  $\frak r\,\colon \Delta_m\to \F^{^{\rm sc}}$ 
are the pairs $(\nu,\delta)$ formed by an {\it order-preserving map\/} 
or, equivalently, a functor  $\delta\,\colon \Delta_m\to \Delta_n$ and by a natural isomorphism 
$\nu\,\colon \frak q\circ\delta\cong \frak r\,,$ the composition being defined by the composition 
of maps and of natural isomorphisms [6,~A2.8]. 

\medskip
£A7.  We say that an {\it $\F^{^{\rm sc}}\-$chain\/} $\frak q\,\colon \Delta_n\to \F^{^{\rm sc}}$ is {\it fully
normalized\/} in~$\F$ if~$\frak q (n)$ is fully normalized in $\F$ and if, moreover, setting
$P' = N_P \big(\frak q (n)\big)$ and~$\F' = N_\F \big(\frak q (n)\big)$, whenever $n\ge 1$ the $\F'\-$chain
$\frak q'\,\colon \Delta_{n-1}\to \F'$ mapping $i\in \Delta_{n-1}$ on the {\it image\/} of $\frak q
(i\bullet n)\,,$ and the $\Delta_{n-1}\-$morphisms on the cor-responding inclusion maps,
is {\it fully~normalized\/} in~$\F'$  [6,~2.18]; note that, by [6, Proposition~2.7], any $\F^{^{\rm sc}}\-$chain
admits a $\ch^* (\F^{^{\rm sc}})\-$isomorphic $\F\-$chain {\it fully normalized\/} in $\F\,.$
Moreover, if $\frak q$ is fully normalized in $\F$ and  $n\ge 1\,,$ we inductively define [6,~2.19]
$$N_P (\frak q) = N_{P'}(\frak q')\qq  N_\F (\frak q) = N_{\F'}(\frak q')
\eqno £A7.1,$$
and it follows from [6, Proposition~2.16] that~$N_\F (\frak q)$ is a Frobenius $N_P (\frak q)\-$ca-tegory;
actually, according to [6, Lemma~2.17] and denoting by $\F (\frak q)$ the image in $\F\big(\frak q (n)\big)$ 
of the group of natural automorphisms of $\frak q\,,$ $\frak q (n)$ is {\it fully $\F (\frak q)\-$normalized\/} in~$\F$ 
and we have
$$N_P (\frak q) = N_P^{\F (\frak q)}\big(\frak q (n)\big)\qq  N_\F (\frak q) =
N_{\F}^{\F (\frak q)}\big(\frak q (n)\big)
\eqno £A7.2.$$
Recall that we have a canonical functor  [6, Proposition~A2.10]
$$\aut_{\F^{^{\rm sc}}} : \ch^*(\F^{^{\rm sc}})\too \Gr
\eqno £A7.3\phantom{.}$$
mapping  any $\F^{^{\rm sc}}\-$chain $\frak q\,\colon \Delta_n\to \F^{^{\rm sc}}$
on~$\F (\frak q)\,.$

\medskip
£A8. We define a {\it folded Frobenius category\/} as a triple $(P,\F,\widehat\aut_{\F^{^{\rm sc}}})$ formed 
by a finite $p\-$group  $P\,,$ by a Frobenius $P\-$category $\F$ and   by the choice  of a functor 
$$\widehat\aut_{\F^{^{\rm sc}}} : \ch^*(\F^{^{\rm sc}})\too k^*\-\Gr
\eqno £A8.1\phantom{.}$$
lifting $\aut_{\F^{^{\rm sc}}}\,;$ note that, for any finite $k^*\-$group $\hat G$ and any block $b$ of $\hat G\,,$ 
denoting by $P$ a defect $p\-$subgroup of $b\,,$ Theorem~11.32 in~[6] guarantees de existence
of a {\it folded Frobenius category\/} $(P,\F_{(b,\hat G)},\widehat\aut_{(\F_{(b,\hat G)})^{^{\rm sc}}})\,.$
{\it Mutatis mutandis\/}, we consider the category $\ch^*(\F^{^{\rm rd}})$ and the canonical functor
$$\aut_{\F^{^{\rm rd}}} : \ch^*(\F^{^{\rm rd}})\too \Gr
\eqno £A8.2;$$
then, it follows from Lemma~£A5 above and from the following result [9, Theorem~2.9] that, for any  subgroup $U$ of $P$ 
fully $K\-$normalized in~$\F$ for some subgroup $K$ of~${\rm Aut}(U)\,,$ we still get a
{\it folded Frobenius category\/} $N_{(P,\F,\widehat\aut_{\F^{^{\rm sc}}})}^K(U)$ formed by the $p\-$group
$N_P^K (U)\,,$ the Frobenius $N_P^K (U)\-$category $N_\F^K (U)$ and the unique functor
$$\widehat\aut_{N_\F^K(U)^{^{\rm sc}}} : \ch^*\big(N_\F^K(U)^{^{\rm sc}}\big)\too k^*\-\Gr
\eqno £A8.3\phantom{.}$$
extending the restriction of $\widehat\aut_{\F^{^{\rm sc}}}$ to $ \ch^*\big(N_\F^K(U)^{^{\rm rd}}\big)\,.$

\bigskip
\noindent
{\bf Theorem~£A9.} {\it Any functor $\widehat\aut_{\F^{^{\rm rd}}}$
lifting $\aut_{\F^{^{\rm rd}}}$ to the category $k^*\-\Gr$ can be extended to a unique functor
lifting $\aut_{\F^{^{\rm sc}}}\,.$
$$\widehat\aut_{\F^{^{\rm sc}}} : \ch^* (\F^{^{\rm sc}})\too k^*\-\Gr 
\eqno  £A9.1\phantom{.}$$\/}

\par
£A10. On the other hand, let us call {\it $k^*\-$localizer\/} any $k^*\-$group $\hat L$ with finite
$k^*\-$quotient $L$ fulfilling $C_L \big(\Bbb O_p(L)\big) = Z\big(\Bbb O_p(L)\big)\,;$
note that $1$ is the unique block of $\hat L\,.$ Following Dade, let us call {\it radical chain\/} of $\hat L$
any subset $\frak r$ of $p\-$subgroups of $\hat L$ which is {\it totally ordered\/} by the inclusion and,
 for any $R\in \frak r\,,$ fulfills $R = \Bbb O_p\big(N_{\hat L}(\frak r^R)\big)$ where $\frak r^R$ is the subset of 
 $\frak r$ of all the elements contained in $R\,;$ note that any element of $\frak r$ contains $\Bbb O_p(L)$ and that 
$\frak r$ can be identified with a {\it regular\/} $\F_{\!(1,\hat L)}^{^{\rm sc}}\-$chain. Now, a $\Bbb Q\-$valued function $f$  
defined over the set of {\it isomorphism classes\/} of  folded Frobenius categories is called {\it radical\/} whenever there exists a $\Bbb Q\-$valued function $f^*$ defined over the set of isomorphism classes
 of {\it $k^*\-$localizers\/} such that for any {\it folded Frobenius category\/}  
 $(P,\F,\widehat\aut_{\F^{^{\rm sc}}})$ we have
 $$f(P,\F,\widehat\aut_{\F^{^{\rm sc}}})  = \sum_{R} f^*\big(\hat L_\F(R)\big)
 \eqno £A10.1$$
 \eject
 \noindent
 where $R$ runs over a set of representatives {\it fully normalized\/} in $\F$ for the set of $\F\-$isomorphism 
 classes of $\F\-$radical subgroups of $P$ and, for such an~$R\,,$ $L_\F (R)$ denotes the $\F\-$localizer of~$R$
 [6,~Theorem~18.6] and $\hat L_\F(R)$
 is the {\it $k^*\-$localizer\/}  coming from the {\it pull-back\/} 
 $$\matrix{\widehat\aut_{\F^{^{\rm sc}}} (R)&\too & \F (R)\cr
\uparrow&&\uparrow\cr
\hat L_\F(R)&\too & L_\F(R) \cr}
\eqno £A10.2.$$

\medskip
\noindent
{\bf Theorem~£A11.} {\it A  $\Bbb Q\-$valued function $f$  defined over the set of {\it isomorphism classes\/} of  folded Frobenius categories is radical if and only if for any folded Frobenius category $(P,\F,\widehat\aut_{\F^{^{\rm sc}}})$
we have 
$$f(P,\F,\widehat\aut_{\F^{^{\rm sc}}}) = \sum_{\frak q} (-1)^{\ell (\frak q)} f\big(N_P(\frak q),N_\F(\frak q),
\widehat\aut_{N_\F(\frak q)^{^{\rm sc}}}\big)
\eqno £A11.1\phantom{.}$$
where $\frak q$ runs over a set of representatives, fully normalized in $\F\,,$ for the set of $\F\-$isomorphism classes of
regular $\F^{^{\rm sc}}\-$chains. In this case, for any $k^*\-$localizer $\hat L\,,$ choosing a Sylow $p\-$subgroup $Q$
of $\hat L$ we have
$$f^*(\hat L) = \sum_\frak r (-1)^{\ell(\frak r)} f(N_Q(\frak r),\F_{(1,N_{\hat L}(\frak r))},
\widehat\aut_{\F_{(1,N_{\hat L}(\frak r))}^{^{\rm sc}}})
\eqno £A11.2\phantom{.}$$
where $\frak r$ runs over a set of representatives, contained in $Q$ and fully normalized in $\F_{\!(1,\hat L)}\,,$ 
for the set of  $\hat L\-$conjugacy classes of radical chains of $\hat L$ such that~$\frak r (0) = \Bbb O_p(\hat L)\,.$\/}

\medskip
\noindent
{\bf Proof:} Firstly assume that $f$ fulfills all the equalities A11.1; then, we claim that it suffices to choose
the function $f^*$ defined by the equalities A11.2; that is to say, we claim that for any folded Frobenius category $(P,\F,\widehat\aut_{\F^{^{\rm sc}}})$
we have 
$$\eqalign{&f(P,\F,\widehat\aut_{\F^{^{\rm sc}}})\cr
& = \sum_R\,\sum_\frak r (-1)^{\ell(\frak r)} f(N_{N_P (R)}(\frak r),
\F_{(1,N_{\hat L_\F(R)}(\frak r))},\widehat\aut_{\F_{(1,N_{\hat L_\F(R)}(\frak r))}^{^{\rm sc}}})\cr}
\eqno £A11.3\phantom{.}$$
 where $R$ runs over a set of representatives {\it fully normalized\/} in $\F$ for the set of $\F\-$isomorphism 
 classes of $\F\-$radical subgroups of $P$ and, for such an $R\,,$  $\frak r$~runs over a set of representatives,
 contained in $N_P (R)$ and fully normalized in~$\F_{\!(1,\hat L_\F(R))}\,,$ for the set of  $\hat L_\F(R)\-$conjugacy classes of radical chains of $\hat L_\F(R)$ such that $\frak r (0) = R\,.$
 
 \smallskip
 But, it is quite clear that such an $\frak r$ can be considered as a {\it regular\/} $\F^{^{\rm sc}}\-$chain which is also
fully normalized in $\F\,,$ and that two of them $\frak r$ and $\frak r'$ fulfilling $\frak r (0) = \frak r'(0)$ are 
$\hat L_\F\big(\frak r(0)\big)\-$conjugate if and only if they are $\F\-$isomorphic; moreover, we clearly have
[6,~Theorem~18.6]
$$N_{N_P (R)}(\frak r) = N_P(\frak r)\qq \F_{(1,N_{\hat L_\F(R)}(\frak r))} \cong N_\F (\frak r) 
\eqno £A11.4.$$
\eject
\noindent
That is to say, the sum
in the right-hand member of equality~£A11.3 contains all the terms of the sum in the right-hand member 
of equality~£A11.1 and therefore, in order to prove our claim, it suffices to show that the sum of the remaining terms
is equal to zero. Hence, our claim follows from Lemmas~£A13 and~£A14 below since the set of isomorphism
classes of the remaining terms coincides with the set $\frak N_0^{^\F} - \frak R^{^\F}$ defined below.

\smallskip
Conversely, assume that $f$ is radical; first of all, we prove that $f$ determines $f^*\,;$ it suffices to show that
 a $\Bbb Q\-$valued function $g$ defined over the set of isomorphism classes
 of {\it $k^*\-$localizers\/} $\hat L$ vanish if  any {\it folded Frobenius category\/}  
 $(P,\F,\widehat\aut_{\F^{^{\rm sc}}})$ fulfills 
 $$0  = \sum_{R} g\big(\hat L_\F(R)\big)
 \eqno £A11.5$$
where $R$ runs over a set of representatives {\it fully normalized\/} in $\F$ for the set of $\F\-$isomorphism  classes of $\F\-$radical subgroups of $P\,;$ we argue by induction on~$\vert L\vert\,.$
 
 \smallskip
 We choose a Sylow $p\-$subgroup $Q$ of $L$ and apply equality~£A11.5 to the folded Frobenius category
 $(Q,\F_{\!(1,\hat L)},\widehat\aut_{\F_{\!(1,\hat L)}})\,,$ so that we have
 $$0  = \sum_{R} g\big(\hat L_{\F_{\!(1,\hat L)}}(R)\big)
 \eqno £A11.6$$
where $R$ runs over a set of representatives {\it fully normalized\/} in $\F_{\!(1,\hat L)}$ for the set of 
$\F_{\!(1,\hat L)}\-$isomorphism  classes of $\F_{\!(1,\hat L)}\-$radical subgroups of $Q\,;$
but, it follows from Lemma~£A5 that such an $R$ contains $\Bbb O_p(L)\,;$ moreover, it follows from
[6,~Corollary~3.6 and~Theorem~18.6] that we have
$$\hat L_{\F_{\!(1,\hat L)}}(R)\cong N_{\hat L} (R)
\eqno £A11.7\phantom{.}$$
and it is clear that $\Bbb O_p(L)$ is  an  $\F_{\!(1,\hat L)}\-$radical subgroup of 
$Q\,;$ since $R\not= \Bbb O_p(L)$ forces $\vert N_L (R)\vert < \vert L\vert\,,$ 
the induction hypothesis implies that all the terms but one vanish in the right hand member of equality~£A11.6, so that we also obtain~$g(\hat L) = 0\,.$

 \smallskip
 Finally, we have to prove that equality~£A11.1 holds; let  $(P,\F,\widehat\aut_{\F^{^{\rm sc}}})$ 
 be a folded Frobenius category; then, for any $\F^{^{\rm sc}}\-$chain $\frak q\,\colon \Delta_n \to \F^{^{\rm sc}}$
 fully normalized in $\F\,,$ we have the  folded Frobenius category $\big(N_P(\frak q),N_\F(\frak q),
\widehat\aut_{N_\F(\frak q)^{^{\rm sc}}}\big)$ and therefore we still have
$$f\big(N_P(\frak q),N_\F(\frak q), \widehat\aut_{N_\F(\frak q)^{^{\rm sc}}}\big) = 
\sum_{R} f^*\big(\hat L_{N_\F(\frak q)} (R)\big)
\eqno £A11.8\phantom{.}$$
where $R$ runs over a set of representatives {\it fully normalized\/} in $N_\F(\frak q)$ for the set of $N_\F(\frak q)\-$isomorphism  classes of $N_\F(\frak q)\-$radical subgroups of $N_P(\frak q)\,;$\break
\eject
\noindent
 consequently, we get
$$\eqalign{\sum_{\frak q} (-1)^{\ell (\frak q)} f\big(N_P(\frak q),N_\F(\frak q),
&\widehat\aut_{N_\F(\frak q)^{^{\rm sc}}}\big)\cr 
&= \sum_{\frak q} (-1)^{\ell (\frak q)}\sum_{R} f^*\big(\hat L_{N_\F(\frak q)} (R)\big)\cr} 
\eqno £A11.9\phantom{.}$$
where $\frak q$ runs over a set of representatives, fully normalized in $\F\,,$ for the set of $\F\-$isomorphism classes
of {\it regular\/} $\F^{^{\rm sc}}\-$chains and, for such a $\frak q\,,$  $R$~runs over a set of representatives 
{\it fully normalized\/} in $N_\F(\frak q)$ for the set of  $N_\F(\frak q)\-$iso-morphism  classes of 
$N_\F(\frak q)\-$radical subgroups of $N_P(\frak q)\,.$

\smallskip
 Once again, it follows from Lemma~£A13 below that
it suffices to consider the sum whenever $\frak q$ belongs to $\frak N^{^\F}\,;$ moreover, we may assume that
$\frak q(i\!-\!1)$ is contained in $\frak q (i)$ and that $\frak q (i\!-\!1\bullet i)$ is the inclusion map 
for any $1\le i\le n\,;$ in this case, since $\frak q (i)\i N_P (\frak q)$ for any $i\in \Delta_n\,,$ we actually have
$$\frak q (i)\i \Bbb O_p\big(\hat L_\F (\frak q)\big)
\eqno £A11.10$$
 for any $i\in \Delta_n\,;$ but, it follows from Lemma~£A5 that $R$ contains  $\Bbb O_p\big(\hat L_\F (\frak q)\big)\,;$
in particular, if $R= \frak q(0)$ then $n = 0\,;$ otherwise, either we have $R\not= \frak q (n)$ and we consider the
{\it regular\/} $\F^{^{\rm sc}}\-$chain $\frak q^\tau\,\colon \Delta_{n+1}\to \F^{^{\rm sc}}$ extending $\frak q$
and mapping $n+1$ on $R$ and the $\Delta_n\-$morphism $n\bullet n+1$ on the corresponding inclusion map, 
or we have $R= \frak q (n)$ 
and we can consider the  restriction $\frak q^\tau\,\colon \Delta_{n-1}\to \F^{^{\rm sc}}$ of~$\frak q\,.$

\smallskip
In both cases, $\frak q^\tau$ belongs to~$\frak N^{^\F}$ and we have $R\not= \frak q^\tau (0)\,;$ moreover, up to replacing
$R$ and $\frak q$ by their image through a suitable $\F^{^{\rm sc}}\-$morphism $R\to P\,,$ we may assume that
$\frak q^\tau$ is fully normalized in $\F$ and then it is easily checked that we get a $k^*\-$isomorphism
$$\hat L_{N_\F(\frak q)} (R)\cong \hat L_{N_\F(\frak q^\tau)} (R)
\eqno £A11.11;$$
 consequently, since we have $(\frak q^\tau)^\tau = \frak q\,,$ in the sum of the right-hand member in~£A11.9 
 only remain the terms where $\frak q (0) = R$ and $n = 0\,,$ and in this case we have  $\hat L_{N_\F(\frak q)} (R)
 = \hat L_\F (R)\,.$ In conclusion, we obtain
 $$\sum_{\frak q} (-1)^{\ell (\frak q)} f\big(N_P(\frak q),N_\F(\frak q),\widehat\aut_{N_\F(\frak q)^{^{\rm sc}}}\big)
 = \sum_R f^*\big(\hat L_\F (R)\big)
\eqno £A11.12\phantom{.}$$
where $\frak q$ and $R$ respectively run over  sets of representatives fully normalized in $\F$ for the 
set of~$\F\-$isomorphism classes of {\it regular\/} $\F^{^{\rm sc}}\-$chains and  for the set of $\F\-$isomorphism 
 classes of $\F\-$radical subgroups of $P\,,$ so that the right-hand member coincides with~$f(P,\F,\widehat\aut_{\F^{^{\rm sc}}})\,.$ We are done.

\medskip
£A12. Let $(P,\F,\widehat\aut_{\F^{^{\rm sc}}})$ be a folded Frobenius category.
For any $i\in \Bbb N\,,$ we denote by $\frak N_i^{^\F}$ the set of $\F\-$isomorphisms classes of {\it regular\/}
$\F^{^{\rm sc}}\-$chains $\frak q\,\colon \Delta_n\to \F^{^{\rm sc}}$ such that, for any $j\in \Delta_n\,,$
if $j < i$ then the image of $\frak q (j)$ in $\frak q (n)$ is {\it normal\/}; moreover,
setting $\frak N^{^\F} = \cap_{i\in \Bbb N}\, \frak N_i^{^\F}\,,$ we denote by $\frak R_i^{^\F}$ the set of 
$\F\-$isomorphisms classes in $\frak N^{^\F}$ of $\F^{^{\rm sc}}\-$chains $\frak q\,\colon \Delta_n\to \F^{^{\rm sc}}$ 
such that, for any $j\in \Delta_n\,,$ if $j < i$ then $\frak q(j)$ is $N_\F (\frak q^j)\-$radical where 
$\frak q^j\,\colon \Delta_j\to \F^{^{\rm sc}}$ is the restriction of $\frak q\,,$ up to replacing
$\frak q^j$ by an $\F\-$isomorphic $\F^{^{\rm sc}}\-$chain fully normalized in $\F\,;$ finally, we set
$\frak R^{^\F} = \cap_{i\in \Bbb N}\, \frak R_i^{^\F}\,.$ It is clear that the stabilizer ${\rm Aut}(P)_\F$ of $\F$
in ${\rm Aut}(P)$ acts on $\frak N_i^{^\F}$ and on $\frak R_i^{^\F}$ for any $i\in \Bbb N\,.$

\bigskip
\noindent
{\bf Lemma~£A13.} {\it With the notation above, for any $i\ge 1$ there is an ${\rm Out}(P)_\F\-$ stable involution $\tau_i$
of the set $\frak N_{i-1}^{^\F} -\frak N_i^{^\F}$ such that, if the isomorphism class $\tilde\frak q$ 
of an $\F^{^{\rm sc}}\-$chain  $\frak q$ fully normalized in $\F$ belongs to $\frak N_{i-1}^{^\F} -\frak N_i^{^\F}\,,$ 
then we have
$$\eqalign{N_{(P,\F,\widehat\aut_{\F^{^{\rm sc}}})} (\frak q)&\cong N_{(P,\F,\widehat\aut_{\F^{^{\rm sc}}})}(\frak  q^{\tau_i})\cr
\hat L_\F (\frak q)\cong \hat L_\F(\frak q^{\tau_i})\quad &and\quad 
\vert \ell (\frak q) - \ell (\frak q^{\tau_i})\vert = 1\cr}
\eqno £A13.1\phantom{.}$$
for a choice in $\tau_i (\tilde\frak q)$ of  an $\F^{^{\rm sc}}\!\-$chain $\frak q^{\tau_i}$ fully normalized in $\F\,.$\/}

\medskip
\noindent
{\bf Proof:} We may assume that $\frak N_{i-1}^{^\F} -\frak N_i^{^\F}\not= \emptyset$ and let 
$\frak q\,\colon \Delta_n\to \F^{^{\rm sc}}$ be an  $\F^{^{\rm sc}}\-$chain  fully normalized in $\F$ 
with its isomorphism class $\tilde\frak q$ in this set; we consider the minimal $j\in \Delta_n$ such that 
$i\le j$ and that the image $Q$ of $\frak q(i -1)$ in~$\frak q(j)$ is not normal; then,  $N_{\frak q (j)}(Q)$ 
is a proper subgroup of $\frak q(j)$ containing the image of $\frak q(j-1)\,.$ If $N_{\frak q (j)}(Q)$ 
coincides with this image, we have $i\not= j -1$ and we consider the $\F^{^{\rm sc}}\-$chain 
$\frak q'\,\colon \Delta_{n-1}\to \F^{^{\rm sc}}$ which coincides with $\frak q$ over $\Delta_{j-2}$ 
and maps $\ell \ge j-1$ on $\frak q(\ell +1)\,;$ otherwise, we consider  the $\F^{^{\rm sc}}\-$chain 
$\frak q'\,\colon \Delta_{n+1}\to \F^{^{\rm sc}}$ which coincides with $\frak q$ over $\Delta_{j-1}$ 
and maps $j$ on $N_{\frak q (j)}(Q)$ and $\ell \ge j+1$ on $\frak q(\ell -1)\,.$ In both cases, 
note that the isomorphism class of~$\frak q'$ still belongs to $\frak N_{i-1}^{^\F} -\frak N_i^{^\F}$ and 
that  $j\in \Delta_n$ is also the minimal element  such that $i\le j$ and that the image
of $\frak q'(i-1)$ in $\frak q'(j)$ is not normal.

\smallskip
Let us replace $\frak q'$ by  an isomorphic $\F^{^{\rm sc}}\-$chain $\frak q^{\tau_i}$ fully normalized in $\F\,;$
in both cases, it is easily checked that such an $\F\-$isomorphism induces the following $\F\-$isomorphism, equivalence of categories and  natural isomorphism
$$N_P (\frak q)\cong N_P (\frak q^{\tau_i})\;,\; N_\F (\frak q)\cong N_\F (\frak q^{\tau_i})
\;\, \hbox{and}\;\,\widehat\aut_{N_\F(\frak q)^{^{\rm sc}}}\cong 
\widehat\aut_{N_\F(\frak q^{\tau_i})^{^{\rm sc}}}
\eqno £A13.2;$$
consequently, according to [6,~Theorem~£18.6] and the {\it pull-back\/}~£A10.2, we get
$$\hat L_\F (\frak q)\cong \hat L_\F(\frak q^{\tau_i})
\eqno £A13.3.$$
Thus, it suffices to define $\tau_i$ as the map sending the isomorphism class of $\frak q$ to the isomorphism
class of $\frak q^{\tau_i}\,.$ We are done.

\bigskip
\noindent
{\bf Lemma~£A14.} {\it  With the notation above, for any $i\ge 1$ there is an ${\rm Out}(P)_\F\-$ stable 
involution $\varpi_i$ of the set $\frak R_{i-1}^{^\F} -\frak R_i^{^\F}$ such that, if the isomorphism class $\tilde\frak q$ 
of an $\F^{^{\rm sc}}\-$chain  $\frak q$ fully normalized in $\F$ belongs to $\frak R_{i-1}^{^\F} -\frak R_i^{^\F}\,,$ 
then we have
$$\eqalign{N_{(P,\F,\widehat\aut_{\F^{^{\rm sc}}})} (\frak q)&\cong 
N_{(P,\F,\widehat\aut_{\F^{^{\rm sc}}})}(\frak  q^{\varpi_i})\cr
\hat L_\F (\frak q)\cong \hat L_\F(\frak q^{\varpi_i})\quad &and\quad 
\vert \ell (\frak q) - \ell (\frak q^{\varpi_i})\vert = 1\cr}\eqno £A14.1\phantom{.}$$
for a choice in $\varpi_i (\tilde\frak q)$ of  an $\F^{^{\rm sc}}\-$chain $\frak q^{\varpi_i}$ fully normalized in $\F\,.$\/}

\medskip
\noindent
{\bf Proof:} We may assume that $\frak R_{i-1}^{^\F} -\frak R_i^{^\F}\not= \emptyset$ and let 
$\frak q\,\colon \Delta_n\to \F^{^{\rm sc}}$ be an  $\F^{^{\rm sc}}\-$chain  fully normalized in $\F$ 
with its isomorphism class $\tilde\frak q$ in this set; that is to say,  since $\frak R_{i-1}^{^\F}\i \frak N^{^\F}\,,$ 
$\frak q(i-1)$ is contained in $N_P(\frak q^{i-1})$ and it is not {\it $N_\F (\frak q^{i-1})\-$radical\/} (cf.~£A7); 
thus,  the structural image  of~$\frak q(i-1)$  in~$L_\F (\frak q^{i-1})$ is a proper subgroup of 
$\Bbb O_p\big(L_\F (\frak q^{i-1})\big)\,,$ and we consider the maximal $j\in \Delta_n$ such that $i-1\le j$ and
that  the structural image $Q$  of~$\frak q(j)$  in~$L_\F (\frak q^j)$ is a proper subgroup of 
$R = \Bbb O_p\big(L_\F (\frak q^j)\big)\,.$

\smallskip
First of all, note that $R$ normalizes the structural image of $\frak q^j$ in $L_\F (\frak q^j)\,;$
moreover, if $j < n$ then the structural image of $\frak q(j+1)$ in $L_\F(q^{j+1})$ coincides
with $\Bbb O_P\big(L_\F (\frak q^{j+1})\big)$ and therefore, since we have [6,~2.13.2 and~Proposition~18.16]
$$L_\F (\frak q^{j+1})\cong N_{L_\F (\frak q^j)}(T)
\eqno £A14.2\phantom{.}$$
where $T$ denotes the structural image of $\frak q(j+1)$ in $L_\F(q^j)\,,$ we still have $N_R (T)\i T\,,$
so that $T$ contains $R\,;$ in conclusion, the structural image of $\frak q(j+1)$~in $L_\F(q^j)$
contains $\Bbb O_p\big(L_\F (\frak q^j)\big)$ which properly contains  the structural~image of~$\frak q(j)\,.$
If $j < n$ and $T = R$ then we consider  the $\F^{^{\rm sc}}\-$chain 
$\frak q'\,\colon \Delta_{n-1}\to \F^{^{\rm sc}}$ which coincides with $\frak q$ over $\Delta_j$ 
and maps $\ell \ge j+1$ on~$\frak q(\ell +1)\,;$  otherwise, we consider  the $\F^{^{\rm sc}}\-$chain 
$\frak q'\,\colon \Delta_{n+1}\to \F^{^{\rm sc}}$ which coincides with $\frak q$ over $\Delta_j$ 
and maps $j +1$ on $\Bbb O_p\big(L_\F (\frak q^j)\big)$  and $\ell \ge j+1$ on $\frak q(\ell -1)\,.$
 In both cases,  note that the isomorphism class of~$\frak q'$ still belongs to $\frak R_{i-1}^{^\F} -\frak R_i^{^\F}$ 
and that  $j\in \Delta_n$ is also the maximal element such that $i-1\le j$ and
that  the structural image of~$\frak q'(j)$  in~$L_\F (\frak q'^j)$ is a proper subgroup of 
$\Bbb O_p\big(L_\F (\frak q'^j)\big)\,.$

\smallskip
Let us replace $\frak q'$ by  an isomorphic $\F^{^{\rm sc}}\-$chain $\frak q^{\varpi_i}$ fully normalized in $\F\,;$
in both cases, it is easily checked that such an $\F\-$isomorphism induces the following  isomorphism of
{\it folded Frobenius categories\/}
$$N_{(P,\F,\widehat\aut_{\F^{^{\rm sc}}})} (\frak q)\cong 
N_{(P,\F,\widehat\aut_{\F^{^{\rm sc}}})}(\frak  q^{\varpi_i})
\eqno £A14.3;$$
consequently, according to [6,~Theorem~£18.6] and the {\it pull-back\/}~£A10.2, we get
$$\hat L_\F (\frak q)\cong \hat L_\F(\frak q^{\varpi_i})
\eqno £A14.4.$$
Thus, it suffices to define $\varpi_i$ as the map sending the isomorphism class of $\frak q$ to the isomorphism
class of $\frak q^{\varpi_i}\,.$ We are done.
\eject

\medskip
£A15. For any $k^*\-$group~$\hat G\,,$ recall that we denote by $\G_k (\hat G)$ the 
{\it scalar extensions\/} from $\Bbb Z$ to~$\O$ of the Grothendieck group  of~the categories 
of finite-dimensional  $k_*\hat G\-$modules; it is well-known that we have a {\it contravariant\/} functor
$$\frak g_k : k^*\-{\Gr}\too \O\-\mod
\eqno  £A15.1\phantom{.}$$
mapping $\hat G$ on $\G_k (\hat G)$ and  any $k^*\-$group homomorphism $\hat\varphi\, \,\colon\hat G\to\hat G'$ on the corresponding {\it restriction\/} map. Then, for any {\it folded Frobenius category\/}  
$(P,\F,\widehat\aut_{\F^{^{\rm sc}}})\,,$ we consider the composed functor
$$\ch^* (\F^{^{\rm sc}})\buildrel \widehat\aut_{\F^{^{\rm sc}}} \over{\hbox to
30pt{\rightarrowfill}}  k^*\-\Gr \buildrel \frak g_k \over{\hbox to 20pt{\rightarrowfill}}  \O\-\mod
\eqno  £A15.2\phantom{.}$$
and we define the (modular) {\it Grothendieck group\/} of $(P,\F,
\widehat\aut_{\F^{^{\rm sc}}})$ as the inverse limit
$$\G_k (P,\F,\widehat\aut_{\F^{^{\rm sc}}}) = \lim_{\longleftarrow}\,
(\frak g_k\circ \widehat\aut_{\F^{^{\rm sc}}})
\eqno  £A15.3;$$
at this point, it follows from [9,~Corollary~8.4] suitably adapted and from Theorem~£A11 above that the 
$\Bbb Z\-$valued function $\rm r$ mapping $(P,\F,\widehat\aut_{\F^{^{\rm sc}}})$ on ${\rm rank}_\O
\big(\G_k (P,\F,\widehat\aut_{\F^{^{\rm sc}}})\big)$ is a {\it radical function\/} and, if the Alperin Conjecture
holds, it is easily checked from [1,~Theorem~3.8] and from  Theorem~£A11 above that $\rm r^*$
maps any {\it $k^*\-$localizer\/} $\hat L$ on the number of blocks of {\it defect zero\/} of
the quotient $\hat L/\Bbb O_p(\hat L)\,.$

\bigskip
\noindent
{\bf References}
\bigskip
\noindent
[1]\phantom{.} Reinhard Kn\"orr and Geoffrey Robinson, {\it Some remarks on a
conjecture of Alperin}, Journal of London Math. Soc. 39(1989), 48-60.
\smallskip\noindent
[2]\phantom{.} Burkhard K\"ulshammer and Llu\'\i s Puig, {\it Extensions of
nilpotent blocks}, Inventiones math., 102(1990), 17-71.
\smallskip\noindent
[3]\phantom{.} Gabriel Navarro and Pham Huu Tiep, {\it ``A reduction theorem for the Alperin weight conjecture''\/},
 Inventiones math., 184(2011), 529-565
 \smallskip\noindent
[4]\phantom{.} Llu\'\i s Puig, {\it Pointed groups and  construction of
modules}, Journal of Algebra, 116(1988), 7-129.
\smallskip\noindent
[5]\phantom{.} Llu\'\i s Puig, {\it ``On the Morita and Rickard
equivalences between Brauer blocks''\/}, Progress in Math., 1999, Birkh\"auser,
Basel.
\smallskip\noindent
[6]\phantom{.} Llu\'\i s Puig, {\it ``Frobenius categories versus Brauer blocks''\/}, Progress in Math. 
274(2009), Birkh\"auser, Basel.
\smallskip\noindent
[7]\phantom{.} Llu\'\i s Puig, {\it Block Source Algebras  in p-Solvable
Groups},  Michigan Math. J. 58(2009), 323-328
\smallskip\noindent
[8]\phantom{.} Llu\'\i s Puig, {\it On the reduction of Alperin's Conjecture to the quasi-simple groups\/}, Journal of Algebra, 
328(2011), 372-398
\smallskip\noindent
[9]\phantom{.} Llu\'\i s Puig, {\it Ordinary Grothendieck groups of a Frobenius $P\-$category\/}, Algebra Colloquium, 18(2011), 1-76
\smallskip\noindent
[10]\phantom{.} Llu\'\i s Puig and Zhou Yuanyang, {\it Glauberman correspondents and extensions of nilpotent block algebras\/}, J. London Math. Soc., Mar 2012; doi: 10.1112/jlms/jdr069 
\smallskip\noindent
[11]\phantom{.} Geoffrey Robinson and Reiner Staszewski, {\it More on Alperin's Conjecture\/}, Ast\'erisque, 181-182(1990), 237-255
\smallskip\noindent
[12]\phantom{.} Jacques Th\'evenaz, {\it Locally determined functions and Alperin's Conjecture\/}, Journal of London Math. Soc., 45(1992), 446-468

\bigskip
\noindent

\end